\documentclass[12pt,a4paper]{amsart}
\usepackage{amsmath,amssymb,amsthm}
\usepackage[a4paper,margin=1in]{geometry}
\usepackage{enumitem}
\usepackage{needspace}
\usepackage[dvipsnames]{xcolor}
\usepackage[colorlinks=true,linkcolor=blue!50!black,citecolor=blue!50!black]{hyperref}
\usepackage{comment}

\newtheorem{theorem}{Theorem}[section]
\newtheorem{lemma}[theorem]{Lemma}
\newtheorem{proposition}[theorem]{Proposition}
\newtheorem{corollary}[theorem]{Corollary}
\newtheorem{mainthm}{Main Theorem}

\theoremstyle{definition}
\newtheorem{definition}[theorem]{Definition}
\newtheorem*{definition*}{Definition}

\newtheorem{remark}[theorem]{Remark}
\newtheorem{question}[theorem]{Question}
\newtheorem{conjecture}[theorem]{Conjecture}

\newcommand{\lk}{\operatorname{lk}}
\newcommand{\skel}{\operatorname{skel}}
\newcommand{\fd}{\operatorname{fd}}
\newcommand{\dist}{\operatorname{dist}}
\newcommand{\diam}{\operatorname{diam}}
\newcommand{\ldiam}{\operatorname{\delta}}
\newcommand{\G}{G}
\newcommand{\Chat}{\widehat{\Delta}}
\newcommand{\M}{\mathcal{M}}

\title[Diameter bounds for complexes close to normal or flag]{Diameter bounds for complexes that are \\close to flag or close to normal}
\author{Bruno Benedetti}
\address[BB]{Department of Mathematics, University of Miami, Coral Gables, FL, USA}
\email{bruno@math.miami.edu}
\author{Giuliamaria Menara}
\address[GM]{Institut de Math\'ematiques de Jussieu-Paris Rive Gauche, Paris, France}
\email{giuliamaria.menara@imj-prg.fr}
\date{November 5, 2026}
\makeatletter
\def\@tocline#1#2#3#4#5#6#7{\relax
  \ifnum #1>\c@tocdepth 
  \else
    \par \addpenalty\@secpenalty\addvspace{#2}%
    \begingroup \hyphenpenalty\@M
    \@ifempty{#4}{%
      \@tempdima\csname r@tocindent\number#1\endcsname\relax
    }{%
      \@tempdima#4\relax
    }%
    \parindent\z@ \leftskip#3\relax \advance\leftskip\@tempdima\relax
    \rightskip\@pnumwidth plus4em \parfillskip-\@pnumwidth
    #5\leavevmode\hskip-\@tempdima #6\nobreak\relax
    \leaders\hbox{$\m@th\mkern 4.5mu\hbox{.}\mkern 4.5mu$}\hfil
    \hbox to\@pnumwidth{\@tocpagenum{#7}}\par
    \nobreak
    \endgroup
  \fi}
\makeatother

\begin{document}

\begin{abstract}
The Hirsch problem asks how large the diameter of the dual graph of a simplicial complex can be, in terms of its number of vertices. For normal complexes, the classical bounds are the Larman--Barnette bound, which is linear for fixed dimension, the Kalai--Kleitman bound, which is quasi-polynomial, and the Adiprasito--Benedetti bound $n-d$, which however holds only for flag complexes. After sharpening the Larman--Barnette and the Kalai--Kleitman bounds, we study how the three bounds deteriorate when flagness or normality is weakened. We give three ways to quantify the distance from flagness, and in each case we bound the diameter in a way that interpolates between the Adiprasito--Benedetti and the Larman--Barnette bounds. We also quantify the distance from normality by a parameter $q$, and obtain a bound that is a polynomial of degree $q$ in $n$. Finally, we conjecture a statement on flag layered families that would imply an interpolation between the Adiprasito--Benedetti and the Kalai--Kleitman bounds.
\end{abstract}

\maketitle

\setcounter{tocdepth}{1}
\tableofcontents

\section{Introduction}

Let $\Delta$ be a pure simplicial complex of dimension $d-1$. We write $\diam(\Delta)$ for the diameter of the dual graph $\G(\Delta)$ of $\Delta$, i.e. for how far apart two facets can be, if we walk from one to the other across faces of dimension $d-2$. The \emph{Hirsch problem}, related to the worst-case running time of the simplex method, is the problem of bounding how large $\diam(\Delta)$ can be in terms of the number $n$ of vertices of $\Delta$.

The question is typically asked for \emph{normal} complexes, i.e. complexes $\Delta$ where for every face $\sigma$, including the empty one, the facets containing $\sigma$ induce a connected subgraph of $\G(\Delta)$. In fact, outside the world of normal complexes, the best general bound for the diameter of strongly connected complexes is elementary: $\diam(\Delta)\le\binom{n}{d-1}$, cf.~Lemma~\ref{lem:ridges}. 

For general normal complexes there are two famous bounds available: the Larman--Barnette bound $2^{d-2}n$ \cite{Lar70,Bar74}, \cite[Thm.~3.5]{LMS17}, 
which is linear in $n$ only when $d$ is fixed, and the Kalai--Kleitman bound $n^{\log_2 d+1}-1$ \cite{KK92}, \cite[Thm.~3.12]{San13}, 
which is quasi-polynomial in $n$ and $d$. A possible linear bound $n-d$, proposed by Hirsch in the Fifties, turned out to be false already for polytopes \cite{San12}, although true for flag normal complexes \cite{AB14} as we will see below. However, all known counterexamples to such a linear bound exceed it only by a small constant factor: to our knowledge, no pure normal complex is known whose dual graph has diameter larger than $1.05 n$, and whether such complexes exist is still open \cite[\S5]{San13}. The so-called ``polynomial Hirsch conjecture'' posits the existence of $a, b \in \mathbb{N}$ such that $n^a \cdot d^b$ is a valid upper bound for all normal complexes \cite{Kal08}. 

Our work begins with a Preliminary Part (Section~\ref{sec:LBKK}) on layered families, whose purpose is to refine the Kalai--Kleitman bound for normal complexes, following the technique of Todd \cite{Tod14} and Sukegawa--Kitahara \cite{SK15}. The new bound is $(n-d+1)^{1+\log_2d}-1$ (Corollary~\ref{cor:normal}). Later in Main Theorem~\ref{thm:B} we also improve the Larman--Barnette bound for normal complexes to $2^{d-2}(n-d)$, again with Todd's idea.

The bulk of the paper is then dedicated to measuring how the classical bounds deteriorate when the hypotheses on $\Delta$ are weakened. We do so in two directions; in part (I) we weaken flagness, and in part (II) we weaken normality.

\bigskip
\noindent{\bfseries Part I: Distance from flagness.}\par\smallskip
In 2014, Adiprasito--Benedetti \cite{AB14} showed that if a $(d-1)$-dimensional simplicial complex $\Delta$ is normal and \emph{flag} (i.e. if every minimal nonface has two elements), then $\diam(\Delta)\le n-d$.
It is therefore natural to ask what survives of this theorem for normal complexes that are ``almost flag''. 
Here we give three answers, corresponding to three different ways of capturing the intuition of a distance from flagness.

\bigskip
\textbf{\textup{(Ia)} A gap in the degrees of the Stanley--Reisner ideal.} For any integer $0\le k\le d-2$, consider the property
\[
(\mathrm{G}_k)\qquad \Delta \text{ has no minimal nonfaces of size } 3,4,\dots,d-k.
\]
In other words, being $G_k$ means that every minimal nonface of $\Delta$ is either a size-2 set, or a large missing face with at least $d-k+1$ vertices. For $k=d-2$ the condition $G_k$ is vacuous, so all $(d-1)$-dimensional complexes satisfy $G_{d-2}$. 
At the opposite extreme, $G_0$ requires every clique with at most $d$ vertices to be a face; the complexes with this property are precisely the \emph{subflag} complexes recently introduced by Benedetti--Pavelka \cite{BP26}. All flag complexes are trivially subflag. The converse is false: all graphs are subflag, but only the triangle-free ones are flag. Let now $\Chat$ be the \emph{flag completion} of $\Delta$, that is, the clique complex of its $1$-skeleton.  
Our first result is the following:

\begin{mainthm}\label{thm:A}
Let $\Delta$ be normal pure $(d-1)$-dimensional on $n$ vertices. Suppose $\Chat$ satisfies the same properties. If  $\Delta$ satisfies $(\mathrm{G}_k)$ for some $k \in \{0, \ldots, d-2\}$, then
\begin{gather*}
\diam(\Delta)\;\le\;(n-d)\cdot \ldiam_k(\Delta),
\end{gather*}
where $\ldiam_0(\Delta):=1$ and $\ldiam_k(\Delta)$ is the maximal diameter of the dual graph of the link of a face with $d-k-1$ vertices (which is a normal $k$-complex). Since for $k \ge 1$
\[ 
\ldiam_k(\Delta)\le \min\Big\{\Big\lfloor\frac{\binom{n-d+k+1}{k}-(k+1)}{k}\Big\rfloor,\;2^{k-1}(n-d),(n-d+1)\tbinom{\lfloor\log_2(n-d+1)\rfloor+k}{k}\Big\}\]
in particular for $k\ge1$, $\diam(\Delta)\le\min\{2^{k-1}(n-d)^2,\,(n-d)(n-d+1)^{1+\log_2(k+1)},\,n^{k+1}\}$. 
\end{mainthm}

For $k=0$ Main Theorem A boils down to the Adiprasito--Benedetti bound: in fact, the assumptions force $\Delta=\Chat$, so $\Delta$ is flag and $\diam(\Delta)\le n-d$.

\bigskip
\textbf{\textup{(Ib)} Flag deficiency.} The next setup investigates the flag property locally, by looking at links. Let us call $\Delta$ \emph{$k$-locally flag} if $\lk_\Delta(\sigma)$ is flag for every face $\sigma$ with $|\sigma|\ge k$. Thus $0$-locally flag means flag, and $1$-locally flag (or simply ``locally flag'') means that all vertex links are flag. Note that every $k$-locally flag complex is also $(k+1)$-locally flag. The \emph{flag deficiency} of $\Delta$ is
\[ \fd(\Delta) \: := \: \min\{k : \Delta\text{ is }k\text{-locally flag}\} = 1+\max\{|\sigma| : \sigma\in \Delta,\ \lk_\Delta(\sigma)\text{ is not flag}\},\]
where by convention $\max\emptyset :=-1$. Hence $\fd(\Delta)=0$ if and only if $\Delta$ is flag, and $\fd(\Delta)\le 1$ if and only if $\Delta$ is locally flag. As one-dimensional links are easy to handle, we also set
\begin{align*}
\fd^{*}(\Delta) := 1+\max\{|\sigma| : \sigma\in \Delta,\ \lk_\Delta(\sigma)\text{ is not flag},\ \dim\lk_\Delta(\sigma)\ge 2\},
\end{align*}
with $\max\emptyset:=-1$. Since $\dim\lk_\Delta(\sigma)=d-1-|\sigma|$, we have $\fd^{*}(\Delta)\le\min\{\fd(\Delta),\,d-2\}$. Moreover, for $d\ge 3$,  $\fd^{*}(\Delta)=0$ if and only if $\Delta$ is flag. Here is our second result:

\begin{mainthm}\label{thm:B}
Let $\Delta$ be normal pure $(d-1)$-dimensional, $d\ge 2$, on $n$ vertices. Then
\[
\diam(\Delta)\le 2^{\fd^{*}(\Delta)}\,(n-d)\le 2^{\min\{\fd(\Delta),\,d-2\}}\,(n-d).
\]
In particular, as $\fd^{*}(\Delta)\le d-2$ always, $\diam(\Delta)\le 2^{d-2}(n-d)$  for every normal complex.\\
Finally, if $\Delta$ is normal and locally flag, then $\diam(\Delta)\le 2(n-d)$.
\end{mainthm}

Main Theorem B interpolates between the Adiprasito--Benedetti bound, which is the case $\fd^{*}(\Delta)=0$, and a slightly sharper form of the Larman--Barnette bound, which corresponds to $\fd^{*}(\Delta)=d-2$.  
The idea of interpolating locally between these bounds goes back to Novik \cite{LMS17}, who first achieved an interpolation result using the ``banner numbers'' from \cite{BjV15,KN13}.
Since \[\fd^{*}(\Delta)\le\fd(\Delta)\le \beta(\Delta)-2,\] where $\beta(\Delta)$ is Klee--Novik's banner number of $\Delta$,  Main Theorem~\ref{thm:B} is stronger than Novik's result $\diam(\Delta)\le 2^{\beta(\Delta)-2}n$ mentioned in \cite[Thm.~3.9]{LMS17}. In Section \ref{sec:examples} we also construct complexes for which Novik's bound is $2^{d-1}n$, while  Main Theorem~\ref{thm:B} gives the better bound $2(n-d)$. 

\bigskip
\textbf{\textup{(Ic)} Almost disjoint minimal nonfaces.} For our third bound we adopt a different viewpoint, based on how minimal nonfaces intersect. We start by proving that if the minimal nonfaces of an arbitrary strongly connected complex $\Delta$  are pairwise disjoint, then $\Delta$ is a join of a simplex with simplex boundaries; thus $\diam(\Delta)=n-d$ exactly. We then perturb this situation by allowing overlaps inside a small, fixed set of $g$ vertices, at the cost of a factor $2^g$ in the diameter bound. Strictly speaking, this approach does not yield a ``distance from flagness'', since flag complexes may very well have missing edges, which may also cover all vertices: think of the one-dimensional cycle $C_5$. However, in the same spirit of Main Theorems \ref{thm:A} and \ref{thm:B} above, we achieve the following bound:

\begin{mainthm}\label{thm:C}
Let $\Delta$ be pure of dimension $d-1$ on vertex set $[n]$, with connected dual graph. Let $g \in \mathbb{N}$. Suppose there is $X\subseteq [n]$ with $|X|=g$ such that $M\cap M'\subseteq X$ for any two distinct minimal nonfaces $M,M'$ of $\Delta$. Then
\[
\diam(\Delta)\;\le\; \kappa\,(n-d+1)-1\;\le\; 2^g(n-d+1)-1,
\]
where $\kappa\le 2^g$ is the number of distinct sets $F\cap X$, with $F$ a facet of $\Delta$.
\end{mainthm}

No normality is needed, and the overlaps $M\cap M'$ may have any size. 
As we said, there are flag complexes with $g$ large, but also   complexes with $g=0$ that look ``far from flag'' from the perspective of Main Theorems \ref{thm:A} and \ref{thm:B}; thus it is clear that Main Theorems \ref{thm:A}, \ref{thm:B}, \ref{thm:C} are logically independent.

\bigskip
\Needspace{6\baselineskip}\medskip\noindent{\bfseries Part II: Distance from normality.}\par\smallskip
So far we studied the distance from flagness; we now wish to study the distance from normality. Since normality asks that the facets containing any face $\sigma$ induce a connected subgraph of $\G(\Delta)$, we are now going to quantify this property by requiring it only for faces with at least $q$ vertices.

\begin{definition*}
Let $1\le q\le d$. A pure $(d-1)$-dimensional complex $\Delta$ is \emph{$q$-normal} if for every face $\sigma$ with $|\sigma|\ge q$ the facets containing $\sigma$ induce a connected subgraph of $\G(\Delta)$. The \emph{normality index} $\nu(\Delta)$ is the least $q$ for which $\Delta$ is $q$-normal.
\end{definition*}

Thus normal means $1$-normal, a $q$-normal complex is $q'$-normal for every $q'\ge q$, and every complex with connected dual graph is $(d-1)$-normal, since the facets containing a face with $d-1$ vertices are pairwise adjacent. Our fourth result is:

\begin{mainthm}\label{thm:D}
Let $\Delta$ be a pure $(d-1)$-dimensional complex on $n$ vertices with connected dual graph, and let $q\ge\nu(\Delta)$ with $1\le q\le d-1$. Put $r=d-q$. Then
\[
\diam(\Delta)\ \le\ 3\,\min\Big\{2^{h}\,\mu_q(\Delta),\ \ B_{r,\,n-d+1}\,\big(\mu_q(\Delta)+1\big)\Big\},
\]
where
\[
\mu_q(\Delta)=\frac{f_{q-1}-\binom dq}{\binom{d-1}{q-1}},\quad B_{r,x}=\binom{\lfloor\log_2x\rfloor+r-1}{r-1},\quad h=\max\{0,\min\{\fd(\Delta)-q,\,r-2\}\}.
\]
\end{mainthm}

For fixed $d$ and $q$, $\mu_q(\Delta)$ is bounded by a polynomial of degree $q$ in $n$: indeed 
\[ n-d\le\mu_q(\Delta)\le\binom nq/\binom{d-1}{q-1}.\] Therefore, the inequality $\diam(\Delta)\le3B_{r,\,n-d+1}(\mu_q(\Delta)+1)$ given by Main Theorem \ref{thm:D} is of Kalai--Kleitman type, while the other inequality,
$\diam(\Delta)\le3\cdot2^h\mu_q(\Delta)$ is of Larman--Barnette type. Note that if $\Delta$ is ``closer to being flag than to being normal'', that is, if $\fd(\Delta)\le \nu (\Delta) \le q$, or more generally if we choose $q = \max \{\fd(\Delta), \nu (\Delta)\}$, then $h=0$. So the first estimate has no exponential factor. In particular Main Theorem \ref{thm:D} specializes to a nice upper bound for the diameter of strongly connected flag complexes, namely,  
\[ \diam(\Delta) \le \frac{3d}{q}\left(\left(\frac nd\right)^{q}-1\right) \quad \textrm{(Corollary \ref{cor:flagq}).}\]

Recall that for flag normal complexes Adiprasito--Benedetti obtained the linear bound $n-d$. In Section~\ref{sec:conj} we conjecture (Conjecture~\ref{conj:flag}) an extension of their theorem to all layerings of flag families, and show that it would imply the following interpolation between Adiprasito--Benedetti and Kalai--Kleitman:

\begin{mainthm}[conditional]\label{thm:E}
Assume Conjecture~{\em \ref{conj:flag}}. Let $\Delta$ be a pure $(d-1)$-dimensional complex on $n$ vertices with connected dual graph and $\fd(\Delta)\le p$.
\begin{enumerate}[label=\textup{(\alph*)}]
\item If $\Delta$ is normal, then $\diam(\Delta)<(n-d+1)\binom{\lfloor\log_2(n-d+1)\rfloor+p}{p}$.
\item If $\Delta$ is $q$-normal, $1\le q\le d-1$ (for instance $q=\nu(\Delta)$), then \\$\diam(\Delta)<3\,(\mu_q(\Delta)+1)\dbinom{\lfloor\log_2(n-d+1)\rfloor+p'}{p'}$, where $p'=\max\{p-q,0\}$.
\end{enumerate}
\end{mainthm}

For $p=0$, part (a) is the Adiprasito--Benedetti bound; for $p=d-1$, which holds for every complex, it is the Kalai--Kleitman bound of Corollary~\ref{cor:normal}, at most $(n-d+1)^{\log_2d+1}$, and no conjecture is needed.

\medskip
\textbf{Organization.} Section~\ref{sec:LBKK} proves the sharper bounds for layered families, Sections~\ref{sec:proofA}--\ref{sec:proofC} prove Main Theorems~\ref{thm:A}--\ref{thm:C}, Sections~\ref{sec:normal}--\ref{sec:conj} prove Main Theorems~\ref{thm:D}--\ref{thm:E}, Section~\ref{sec:examples} has examples, and Section~\ref{sec:questions} lists some open questions.

\newpage

\Needspace{6\baselineskip}\section{Layered families and a sharper version of Kalai--Kleitman}\label{sec:LBKK}

In this section we extend Santos' treatment of layered families \cite[Section~3]{San13} and, following Todd's idea \cite{Tod14}, we express the bounds in terms of $n-d$ instead of $n$. We need this sharper form later in the paper (Main Theorems~\ref{thm:A}, \ref{thm:D} and~\ref{thm:E}), where the bounds are applied to the links of faces, which have few vertices and small dimension.

Here, a \emph{$1$-layered family of $r$-sets} is a set of $r$-subsets of a finite set, partitioned into nonempty layers $L_0,\dots,L_N$, such that for every nonempty $S$ the set of $t$ with $S\subseteq F$ for some $F\in L_t$ is an interval. Unlike the multicomplexes of \cite{San13}, these are families of \emph{sets}. For $m\ge r\ge1$ let $F(m,r)$ be the maximal number $N+1$ of layers of such a family using at most $m$ vertices. Put $F(m,r)=0$ for $m<r$. Santos \cite[Thms.~3.12 and~3.14]{San13} proves the Kalai--Kleitman and Larman--Barnette bounds for connected layer multicomplexes. The proofs below follow his, with two changes: we count layers rather than the length, and we exploit that the sets have $r$ distinct elements, which allows us to replace $m$ by $m-r$ or $m-r+1$. This is the same idea that Todd \cite{Tod14} and Sukegawa and Kitahara \cite{SK15} used to improve the Kalai--Kleitman bound for polytopes.

\begin{theorem}\label{thm:LBKK}
For all $r\ge1$ and $m\ge r$,
\begin{align*}
F(m,r)&\ \le\ \lambda_r(m):=2^{r-1}(m-r)+1,\\
F(m,r)&\ \le\ \kappa_r(m):=(m-r+1)\binom{\lfloor\log_2(m-r+1)\rfloor+r-1}{r-1}.
\end{align*}
Moreover $\kappa_r(m)\le(m-r+1)^{1+\log_2r}$.
\end{theorem}

For comparison, Santos' Theorem~3.14 \cite{San13} gives at most $2^{r-1}(m-1)+1$ layers, and his Theorem~3.12 gives at most $m^{1+\log_2r}$ layers. For $r=1$ all of these equal $m$, which is the exact value.

\begin{proof}[Proof of the first bound]
We use induction on $r$. For $r=1$ the sets are distinct singletons, so $F(m,1)=m=\lambda_1(m)$. Let $r\ge2$ and $m\ge r$, and let $L_0,\dots,L_N$ be the layers of a $1$-layered family of $r$-sets on a set $V$ with $|V|=m$. Let $s_1=0$. If $s_i\le N$, let $l_i\ge s_i$ be the largest index such that $L_{l_i}$ and $L_{s_i}$ share a vertex, and put $s_{i+1}=l_i+1$. The $i$-th block consists of the layers $L_t$ with $s_i\le t\le l_i$; let $k$ be the number of blocks. Choose a vertex $v_i$ common to $L_{s_i}$ and $L_{l_i}$; it lies in all the layers of the block, because the layers containing $v_i$ form an interval. Let $n_i$ be the number of vertices used by the $i$-th block. The sets $F\setminus v_i$, with $F$ in the $i$-th block and $v_i\in F$, form a $1$-layered family of $(r-1)$-sets on $n_i-1$ vertices, in which each layer of the block is represented, so the block has at most $2^{r-2}(n_i-r)+1$ layers by induction, since $n_i-1\ge r-1$ and $\lambda_{r-1}$ is increasing.

A vertex $w$ lies in at most two blocks. Indeed, if $w$ lies in blocks $i$ and $j$ with $j\ge i+2$, then by the interval property $w$ lies in $L_{s_{i+1}}$ and in a layer $L_t$ with $t\ge s_j>l_{i+1}$, contradicting the choice of $l_{i+1}$. The same argument shows that the $a\ge r$ vertices used by $L_0$ lie in the first block only, since otherwise they would contradict the choice of $l_1$. Hence $\sum_in_i\le a+2(m-a)\le2m-r$, and
\begin{align*}
N+1&\ \le\ \sum_{i=1}^k\big(2^{r-2}(n_i-r)+1\big)\ \le\ 2^{r-2}\big(2m-r-kr\big)+k\\
&\ =\ 2^{r-1}(m-r)+1-(k-1)\big(2^{r-2}r-1\big)\ \le\ \lambda_r(m),
\end{align*}
because $2^{r-2}r\ge2$.
\end{proof}

\begin{proof}[Proof of the second bound]
We first show that, for $m\ge r\ge2$,
\begin{equation}\label{eq:rec}
F(m,r)\ \le\ F(\lfloor\tfrac{m-1}2\rfloor,r)+F(\lceil\tfrac{m-1}2\rceil,r)+F(m-1,r-1).
\end{equation}
This is Santos' Lemma~3.11 \cite{San13} in terms of layers. Take a family with layers $L_0,\dots,L_N$ on $m$ vertices. Let $a=\lfloor\frac{m-1}2\rfloor$ and $b=\lceil\frac{m-1}2\rceil$. Let $i$ be the largest integer such that $L_0,\dots,L_{i-1}$ use at most $a$ vertices in total, and $j$ the largest integer such that the last $j$ layers use at most $b$ vertices in total. Then $i\le F(a,r)$ and $j\le F(b,r)$. If $i+j\ge N+1$ we are done. Otherwise, let $k=N+1-i-j\ge1$. The layers $L_0,\dots,L_i$ use at least $a+1$ vertices, and $L_{N-j},\dots,L_N$ use at least $b+1$; since $(a+1)+(b+1)=m+1>m$, some vertex $v$ lies in a layer of index $\le i$ and in a layer of index $\ge N-j$. By the interval property, $v$ lies in each of the $k$ layers $L_i,\dots,L_{N-j}$. Taking links of $v$ in these layers gives a $1$-layered family of $(r-1)$-sets on at most $m-1$ vertices with $k$ layers, so $k\le F(m-1,r-1)$. This proves~\eqref{eq:rec}.

We prove $F(m,r)\le\kappa_r(m)$ by induction on $r$ and then on $m$. For $r=1$, $F(m,1)=m=\kappa_1(m)$. For $m=r$ there is only one $r$-set, and $\kappa_r(r)=1$. Let $r\ge2$ and $m>r$, put $\mu=m-r+1\ge2$ and $L=\lfloor\log_2\mu\rfloor\ge1$. By induction, since $(m-1)-(r-1)+1=\mu$,
\[
F(m-1,r-1)\ \le\ \mu\binom{L+r-2}{r-2}.
\]
For $x\in\{a,b\}$ put $\mu_x=x-r+1$. If $\mu_x\le0$ then $F(x,r)=0$. If $\mu_x\ge1$ then $\mu_x\le\mu/2$, because $x\le m/2$ and $r\ge1$. Hence $\lfloor\log_2\mu_x\rfloor\le L-1$, and by induction $F(x,r)\le\mu_x\binom{L+r-2}{r-1}$. In all cases $F(x,r)\le\frac\mu2\binom{L+r-2}{r-1}$. Adding up with \eqref{eq:rec} and using Pascal's rule,
\[
F(m,r)\ \le\ \mu\left[\binom{L+r-2}{r-1}+\binom{L+r-2}{r-2}\right]=\mu\binom{L+r-1}{r-1}=\kappa_r(m).
\]
Finally $\binom{L+r-1}{r-1}=\prod_{j=1}^L\frac{r-1+j}j\le r^L\le r^{\log_2\mu}=\mu^{\log_2r}$.
\end{proof}

\begin{corollary}[Improved Kalai-Kleitman]\label{cor:normal}
Let $\Delta$ be a pure normal $(d-1)$-dimensional complex on $n$ vertices. Then
\[
\diam(\Delta)\ \le\ (n-d+1)\binom{\lfloor\log_2(n-d+1)\rfloor+d-1}{d-1}-1.
\]
\end{corollary}

\begin{proof}
Let $A$ be an endpoint of a longest shortest path, let $e=\diam(\Delta)$, and let $L_t$ be the set of facets at distance $t$ from $A$, for $0\le t\le e$. For a face $S$, the facets containing $S$ induce a connected subgraph of $\G(\Delta)$, because $\Delta$ is normal, and distances from $A$ change by at most one along an edge; hence the set of $t$ with $S\subseteq F$ for some $F\in L_t$ is an interval. So $L_0,\dots,L_e$ form a $1$-layered family of $d$-sets on $n$ vertices with $e+1$ layers, and the second bound of Theorem~\ref{thm:LBKK} with $r=d$ and $m=n$ applies. (The first bound gives only $2^{d-1}(n-d)$, which is worse than the bound $2^{d-2}(n-d)$ of Main Theorem~\ref{thm:B}.)
\end{proof}

\begin{remark}\label{rem:SK}
Todd \cite{Tod14} proved that a $d$-dimensional polyhedron with $n$ facets has diameter at most $(n-d)^{\log_2d}$, and Sukegawa and Kitahara \cite{SK15} improved this to $(n-d)^{\log_2(d-1)}$. 
These results are stated for polyhedra; as far as we understand, the proofs use steps that are specific to polyhedra, such as the reduction to simple polyhedra, the restriction to the sub-polyhedron defined by some subset of the facets, or Klee's bound $n-3$ for the diameter of $3$-polytopes. We do not know whether these steps extend to normal complexes; as for layered families of sets, a low-dimensional bound corresponding to Klee's fails without flagness assumptions (cf.~Section~\ref{sec:conj}). Our bound $\kappa_r$ above follows the idea of Todd and Sukegawa--Kitahara; it is obtained by the same shift of $m$ by $r$, carried out directly in the recursion for layered families, where the base cases $r=1$ and $m=r$ are trivial; but strictly speaking, it is not a direct consequence of the results in \cite{Tod14} and \cite{SK15}.
\end{remark}

\Needspace{6\baselineskip}\section{Distance from flagness: proof of Main Theorem~\ref{thm:A}}\label{sec:proofA}

\textbf{Notation}. Throughout, $\Delta$ is a (finite) pure simplicial complex of dimension $d-1$ on a vertex set $V$ with $|V|=n$. The \emph{dual graph} $\G(\Delta)$ has the facets of $\Delta$ as nodes, and its edges join facets that are adjacent, i.e. that share exactly $d-1$ vertices; $\diam(\Delta)$ is its diameter. The \emph{link} of a face $\sigma$ is $\lk_\Delta(\sigma)=\{\tau : \tau\cap\sigma=\emptyset,\ \tau\cup\sigma\in \Delta\}$. The \emph{star} of a face $\sigma$ is the smallest subcomplex of $\Delta$ containing all facets that contain $\sigma$. We use freely that $\lk_{\lk_\Delta(\sigma)}(\tau)=\lk_\Delta(\sigma\cup\tau)$, and that $\lk_\Delta(\sigma)$ has the dual graph of the star of $\sigma$, the subgraph of $\G(\Delta)$ induced by the facets containing $\sigma$. Links of pure normal complexes are pure and normal.
A \emph{minimal nonface} is a set $M\subseteq V$ with $M\notin \Delta$ but all proper subsets in $\Delta$; $V$ is always the vertex set of $\Delta$, so minimal nonfaces have at least two elements. For a finite set $Z$, $\langle Z \rangle$ is the full simplex on $Z$ and $\partial\langle Z \rangle$ its boundary.
Throughout this section fix $k$ with $0\le k\le d-2$ and put $s:=d-k-1\ge1$. Then $(\mathrm{G}_k)$ says that $\Delta$ has no minimal nonfaces of size $3,\dots,s+1$, and the links of faces with $s$ vertices are normal of dimension $(d-1)-s=k$. For $k=0$ these links are $0$-dimensional, so their dual graphs are complete.
Let us start with two easy lemmas:

\begin{lemma}\label{lem:gap}
For $1\le s\le d-1$, the following are equivalent: \em
\begin{enumerate}[nosep,label=(\roman*)]
\item {\it $\Delta$ has no minimal nonfaces of size $3,\dots,s+1$};
\item {\it every clique of the $1$-skeleton with at most $s+1$ vertices is a face}; 
\item $\skel_{s}\Delta=\skel_{s}\Chat$.
\end{enumerate}
\end{lemma}

\begin{proof}
If $S$ is a clique with $|S|\le s+1$ that is not a face, a minimal nonface $M\subseteq S$ is a clique, so $|M|\ge 3$, and $|M|\le s+1$. Conversely, a minimal nonface with at least three elements is a clique. The equivalence of \textup{(ii)} and \textup{(iii)} follows from the definition of $\Chat$, since $\skel_s$ consists of the faces with at most $s+1$ vertices.
\end{proof}

\begin{lemma} \label{lem:ridges}
Let $L$ be pure of dimension $t\ge 1$ on $n$ vertices with connected dual graph. Then $\diam(L)\le\big\lfloor(\binom{n}{t}-(t+1))/t\big\rfloor\le\binom{n}{t}$. If $t=1$, then $\diam(L)\le n-2$.
\end{lemma}

\begin{proof}
Along a shortest path $Q_0,\dots,Q_a$ the ridges $Q_{i-1}\cap Q_i$ are distinct: if $Q_{i-1}\cap Q_i=Q_{j-1}\cap Q_j$ with $i<j$, then $Q_{i-1}$ and $Q_j$ are equal or adjacent, and the path could be shortened. A ridge has $t$ vertices, and the same argument shows that a ridge of $Q_i$ is a ridge of $Q_j$ only if $|i-j|\le1$. Hence $Q_0$ contributes $t+1$ ridges and each further $Q_i$ contributes $t$ new ones, so $(t+1)+at\le\binom{n}{t}$. For $t=1$, $L$ is a connected graph and $\G(L)$ its line graph. If edges $e,f$ are disjoint, take endpoints $u\in e$, $w\in f$ at minimal distance $\ell$ in $L$; the $\ell+1$ vertices of a shortest $u$--$w$ path and the other endpoints of $e,f$ are distinct, so $\ell\le n-3$, and $\dist(e,f)\le \ell+1\le n-2$.
\end{proof}

For a pure normal $(d-1)$-complex $\Delta$ define its \emph{$k$-local diameter} by
\[
\ldiam_k(\Delta):=\max\{\diam(\lk_\Delta\sigma) : \sigma\in \Delta,\ |\sigma|=d-k-1\}.
\]
The link of such a $\sigma$ is normal of dimension $k$ on at most $n-d+k+1$ vertices, so Lemma~\ref{lem:ridges} and Main Theorem~\ref{thm:B} (applied to the link, with $d$ replaced by $k+1$; its proof does not use Main Theorem~\ref{thm:A}) give
\begin{equation}\label{eq:ldiam}
\begin{split}\ldiam_k(\Delta)\le \min\Big\{&\Big\lfloor\frac{\binom{n-d+k+1}{k}-(k+1)}{k}\Big\rfloor,\ 2^{k-1}(n-d),\\&(n-d+1)\binom{\lfloor\log_2(n-d+1)\rfloor+k}{k}\Big\},\end{split}
\end{equation}
since a normal $k$-dimensional complex on $m$ vertices has diameter at most $2^{k-1}(m-k-1)$ by Main Theorem~\ref{thm:B}, and at most $(m-k)\binom{\lfloor\log_2(m-k)\rfloor+k}{k}-1$ by Corollary~\ref{cor:normal}, where $m-k\le n-d+1$. (These bounds are for $k\ge1$. For $k=0$ the links of faces with $d-1$ vertices are finite sets of points, and $\ldiam_0(\Delta)=1$.)

\begin{proposition}\label{prop:compare}
Let $\Delta\subseteq D$ be pure $(d-1)$-complexes such that $\skel_{s}\Delta=\skel_{s}D$, with $s=d-k-1\ge1$. Assume $\Delta$ is normal and $\G(D)$ is connected. Then for any facets $A,B$ of $\Delta$,
\[
\dist_{\G(D)}(A,B)\le\dist_{\G(\Delta)}(A,B)\le\ldiam_k(\Delta)\,\dist_{\G(D)}(A,B).
\]
In particular $\diam(\Delta)\le\ldiam_k(\Delta)\,\diam(D)$.
\end{proposition}

\begin{proof}
Facets of $\Delta$ are $(d-1)$-faces of $D$, hence facets of $D$, and $\G(\Delta)$ is an induced subgraph of $\G(D)$; this gives the first inequality.

For the second inequality, we will follow a path in $D$ by carrying an $(s+1)$-vertex face inside each of its facets.
These smaller faces belong to $\Delta$ because the two complexes have the same $s$-skeleton. We choose them so that consecutive ones share at least $s$ vertices. After extending them to facets of $\Delta$, we can replace each step by a path in the link of a common $s$-vertex face, at a cost of at most $\ldiam_k(\Delta)$.
Let $A=F_0,F_1,\dots,F_\ell=B$ be a shortest path in $\G(D)$, with $F_{i+1}=(F_i\setminus a_i)\cup b_i$.
Choose $S_0\subseteq F_0$ with $|S_0|=s+1$ and put $S_{i+1}=S_i$ if $a_i\notin S_i$ and $S_{i+1}=(S_i\setminus a_i)\cup b_i$ otherwise. Then $S_i\subseteq F_i$, $|S_i|=s+1$, $|S_i\cap S_{i+1}|\ge s$, and $S_i\in\skel_{s}D=\skel_{s}\Delta$. 
By purity choose facets $A_i\supseteq S_i$ of $\Delta$, with $A_0=A$ and $A_\ell=B$ (possible since $S_0\subseteq A$ and $S_\ell\subseteq B$). Then $|A_i\cap A_{i+1}|\ge|S_i\cap S_{i+1}|\ge s$.
Pick $\sigma_i\subseteq S_i\cap S_{i+1}$ with $|\sigma_i|=s$. Then $A_i$ and $A_{i+1}$ are joined in $\G(\Delta)$ by a path of length at most $\ldiam_k(\Delta)$ through facets containing $\sigma_i$. Concatenating gives a walk of length at most $\ell\,\ldiam_k(\Delta)$ from $A$ to $B$ in $\G(\Delta)$.
\end{proof}

\begin{proof}[\bf Proof of Main Theorem~\ref{thm:A}]
By Lemma~\ref{lem:gap}, Proposition~\ref{prop:compare} applies with $D=\Chat$. By the Adiprasito--Benedetti bound, $\diam(\Chat)\le n-d$. Combine this with Equation \eqref{eq:ldiam}, and for $k=0$ use $\ldiam_0=1$. For the three simplified forms ($k\ge1$), the first is immediate; for the second, $\binom{L+k}{k}\le(k+1)^L\le (n-d+1)^{\log_2(k+1)}$ for $L=\lfloor\log_2(n-d+1)\rfloor$; for the third note that $n-d+k+1\le n$ because $k\le d-1$, and $\binom{m}{k}\le m^{k}$.
\end{proof}

\begin{remark}\label{rem:A}
The normality of $\Delta$ is used in the proof only for faces with $d-k-1$ vertices (their stars are connected, and their links are then normal as in the proof of Main Theorem~\ref{thm:D}), so Main Theorem~\ref{thm:A} holds as stated if $\Delta$ is only $(d-k-1)$-normal, that is, if $\nu(\Delta)\le d-k-1$. The hypothesis that $\Chat$ be normal is used only to invoke the Adiprasito--Benedetti theorem; indeed, any bound $\diam(\Chat)\le B$ would yield $\diam(\Delta)\le\ldiam_k(\Delta)\cdot B$. For $k\le1$ the normality of $\Chat$ follows from the other hypotheses (Proposition~\ref{prop:normalhat} below); for $k\ge2$ we do not know whether it does (Question~\ref{q:normalhat}).
\end{remark}

\begin{proposition}\label{prop:normalhat}
Let $\Delta$ be pure and normal of dimension $d-1$, satisfying $(\mathrm G_1)$, and assume that $\Chat$ is pure of dimension $d-1$. Then $\Chat$ is normal.
\end{proposition}

\begin{proof}
By Lemma~\ref{lem:gap}, $\Delta$ and $\Chat$ have the same $(d-2)$-skeleton. Let $\sigma\in\Chat$ be a face with at most $d-1$ vertices (including $\sigma=\emptyset$); then $\sigma\in\Delta$. The facets of $\Delta$ containing $\sigma$ are also facets of $\Chat$, they induce a connected subgraph of $\G(\Chat)$ because $\Delta$ is normal, and there is at least one of them because $\Delta$ is pure. Let $F$ be a facet of $\Chat$ containing $\sigma$ that is not a facet of $\Delta$. Choose $v\in F\setminus\sigma$ and set $R=F\setminus v$. Then $R\supseteq\sigma$ has $d-1$ vertices, so $R\in\Delta$, and by purity $R$ lies in a facet $F'$ of $\Delta$. Then $F'\supseteq\sigma$ is adjacent to $F$. Hence all facets of $\Chat$ containing $\sigma$ induce a connected subgraph. 
\end{proof}

\Needspace{6\baselineskip}\section{Distance from flagness: proof of Main Theorem~\ref{thm:B}}\label{sec:proofB}

\subsection{Monotone conservative paths (MCP)} Introduced by Labb\'e, Manneville and Santos \cite{LMS17}, monotone conservative paths (MCP) are a local reformulation of the combinatorial segments of \cite{AB14}. Let $\Delta$ be pure $(d-1)$-dimensional and normal. We use the following facts from \cite{LMS17}: 
\begin{enumerate}[label=\textup{(\roman*)}]
\item Any two facets of $\Delta$ are joined by an MCP \cite[Cor.~2.10, 2.14]{LMS17}.
\item An MCP $\Gamma$ towards a target set $S$ splits into $b$ consecutive blocks $\Gamma_1,\dots,\Gamma_b$ according to its anchors $x_1,\dots,x_b$, where $\Gamma_i$ includes the step in which the anchor changes from $x_i$ to $x_{i+1}$; deleting $x_i$ from $\Gamma_i$ gives an MCP $\Gamma_i'$ in $\lk_\Delta(x_i)$ of the same length \cite[\S2.2, proof of Thm.~3.5]{LMS17}.
\item Say $m$ vertices appear in $\Gamma$ and $m_i$ in $\Gamma_i'$. If  $b=1$ (with the notation above), $m_1=m-1$; if instead  $b\ge2$ then $\sum_i m_i\le 2m-2$. Indeed, the anchor of a facet is its vertex closest to $S$, and the anchors form a shortest vertex path towards $S$, with strictly decreasing distances. A vertex at distance $\ell$ from $S$ that never becomes an anchor can therefore only appear in blocks whose anchor is at distance $\ell$ or $\ell-1$, hence in at most two of the paths $\Gamma_i'$. An anchor $x_i$ is deleted from its own block, and can only appear at the endpoint of the preceding block and in the following one, hence in at most two of the $\Gamma_j'$; for the first and the last anchor, one of these two blocks does not exist. So $\sum_i m_i\le 2(m-2)+2$ when $b\ge2$.
\item In a flag normal complex, MCPs are non-revisiting \cite[Cor.~3.2]{LMS17}, \cite{AB14}: every step introduces a vertex that has not appeared before. Hence an MCP supported on $m$ vertices in a flag normal $(r-1)$-complex has length $\le m-r$.
\end{enumerate}

\begin{definition}
Let $\mathcal H$ be a class of normal complexes that is closed under taking links and contains every normal complex of dimension at most $0$ (including $\{\emptyset\}$). 
The \emph{$\mathcal H$-deficiency} of a normal complex $\Delta$ is $\operatorname{fd}_{\mathcal H}(\Delta):=1+\max\{|\sigma|:\sigma\in \Delta,\ \lk_\Delta\sigma\notin\mathcal H\}$, with $\max\emptyset=-1$. 
Thus $\operatorname{fd}_{\mathcal H}(\Delta)=0$ if and only if $\Delta\in\mathcal H$, and since $\lk_{\lk_\Delta x}(\sigma)=\lk_\Delta(\sigma\cup x)$, 
\begin{equation}\label{eq:depthdrop}
\operatorname{fd}_{\mathcal H}(\lk_\Delta x)\le\max\{0,\operatorname{fd}_{\mathcal H}(\Delta)-1\}\quad\text{for every vertex }x.
\end{equation}
\end{definition}
\begin{theorem}\label{thm:depth}
Suppose there is a constant $A>0$ such that every MCP in any complex of $\mathcal H$ of dimension $r-1$, supported on $m$ vertices, has length at most $A(m-r)$. Then every MCP $\Gamma$ in a normal complex $\Delta$ of dimension $r-1$, if supported on $m$ vertices, satisfies $|\Gamma|\le A\,2^{\operatorname{fd}_{\mathcal H}(\Delta)}(m-r)$.
\end{theorem}

\begin{proof}
We proceed by induction on $p=\operatorname{fd}_{\mathcal H}(\Delta)$; the case $p=0$ is the hypothesis. Let $p\ge1$. Then $\Delta\notin\mathcal H$, so $\dim\Delta\ge1$ and $r\ge2$. Split $\Gamma$ into $b$ blocks as in item (ii) above. Each $\Gamma_i'$ is an MCP in the normal complex $\lk_\Delta(x_i)$ of dimension $r-2$, whose deficiency is at most $p-1$ by \eqref{eq:depthdrop}. By (iii), $\sum_i\big(m_i-(r-1)\big)=m-r$ if $b=1$, and if $b\ge2$,
\[
\sum_i\big(m_i-(r-1)\big)\le 2m-2-b(r-1)\le 2(m-r).
\]
In both cases the sum is at most $2(m-r)$. By induction,
\[
|\Gamma|=\sum_i|\Gamma_i'|\le A\,2^{p-1}\sum_i\big(m_i-(r-1)\big)\le A\,2^{p}(m-r).\qedhere
\]
\end{proof}

\begin{proof}[\bf Proof of Main Theorem~\ref{thm:B}]
Take for $\mathcal H$ the class of normal complexes that are flag or of dimension at most $1$. This class  is closed under taking links, since links of flag complexes are flag and links of complexes of dimension at most $1$ have dimension at most $0$. Any normal complex of dimension $1$ is a connected graph, and its MCPs are shortest paths \cite[proof of Thm.~3.9]{LMS17}, which introduce a new vertex at each step, so their length is at most $m-2$. In dimension zero, a path has length at most $1\le m-1$. Together with item  (iv) above, this shows that MCPs in $\mathcal H$ of dimension $r-1$ supported on $m$ vertices have length at most $m-r$, that is, $A=1$. Moreover $\operatorname{fd}_{\mathcal H}=\fd^{*}$, where $\fd^{*}$ is as defined in the introduction. Also, Theorem~\ref{thm:depth} together with (i), applied with $m\le n$ and $r=d$, gives $\diam(\Delta)\le 2^{\fd^{*}(\Delta)}(n-d)$.
\end{proof}

\begin{remark}\label{rem:nonnormal}
Normality cannot be dropped from Main Theorem~\ref{thm:B}, even when $\fd(\Delta)=1$. Let $\Delta$ be the pure $2$-dimensional complex on $\{0,\dots,6\}$ (so $n=7$ and $d=3$) with the ten facets
\[012,\ 026,\ 056,\ 045,\ 034,\ 234,\ 235,\ 135,\ 136,\ 146.\]
Consecutive facets in this list share an edge, and no other two facets do, so $\G(\Delta)$ is a path with ten vertices and $\diam(\Delta)=9$. Every vertex link is a disjoint union of paths, hence flag, so $\fd(\Delta)\le1$ (in fact $\fd(\Delta)=\fd^{*}(\Delta)=1$, since $013$ is a missing triangle). But $\diam(\Delta)=9>8=2(n-d)$. Of course $\Delta$ is not normal: the facets containing the vertex $1$ are $012,135,136,146$, and the link of $1$ has the isolated edge $02$, so its dual graph is disconnected.
\end{remark}

\subsection{Comparison with banner complexes}
Recall that a clique $T$ of the $1$-skeleton is \emph{critical} if $T\setminus v\in \Delta$ for some $v\in T$. A simplicial complex is \emph{$b$-banner} if every critical clique with at least $b+1$ vertices is a face \cite{KN13}; the \emph{banner number} $\beta(\Delta)$ is the least $b\ge2$ for which $\Delta$ is $b$-banner. Flag is the same as $2$-banner, every $(d-1)$-complex is $(d+1)$-banner, and if $\Delta$ is $b$-banner then each vertex link is $(b-1)$-banner \cite[Lemma~3.8]{LMS17}.

\begin{proposition}\label{prop:banner} Under the previous notation,
\[\fd(\Delta)\le\beta(\Delta)-2.\] 
In particular, Main Theorem~\ref{thm:B} sharpens Novik's bound $\diam(\Delta)\le 2^{\beta(\Delta)-2}n$ for normal complexes $\Delta$.
\end{proposition}

\begin{proof}
Iterating, $\lk_\Delta\sigma$ is $(\beta-|\sigma|)$-banner, hence flag once $|\sigma|\ge\beta-2$. So every face with non-flag link has at most $\beta-3$ vertices.
\end{proof}

\Needspace{6\baselineskip}\section{Distance from flagness: proof of Main Theorem~\ref{thm:C}}\label{sec:proofC}

\begin{lemma}\label{lem:join}
Let $W$ be a finite set, $Z_1,\dots,Z_t\subseteq W$ pairwise disjoint nonempty sets, and
\[
K=\{\tau\subseteq W:\ Z_j\not\subseteq\tau\text{ for all }j\}.
\]
Then $K=\partial\langle Z_1\rangle*\dots*\partial\langle Z_t\rangle*\big\langle W\setminus\textstyle\bigcup_jZ_j\big\rangle$, where $\partial\langle Z\rangle=\{\emptyset\}$ if $|Z|=1$. The facets of $K$ are the sets $W\setminus\{z_1,\dots,z_t\}$ with $z_j\in Z_j$; $K$ is pure with facets of size $|W|-t$; $\G(K)\cong K_{|Z_1|}\,\square\,\cdots\,\square\,K_{|Z_t|}$; and $\diam(K)=\#\{j:|Z_j|\ge2\}\le t$.
\end{lemma}

\begin{proof}
Since the $Z_j$ are disjoint, $\tau$ contains no $Z_j$ if and only if $\tau\cap Z_j\subsetneq Z_j$ for every $j$, with no condition on $\tau\setminus\bigcup Z_j$; this is the join. A maximal such $\tau$ omits exactly one element $z_j$ from each $Z_j$. Two facets indexed by $z,z'$ meet in $|W|-t-\#\{j:z_j\ne z'_j\}$ elements, so they are adjacent if and only if $z,z'$ differ in exactly one coordinate. The diameter of a Cartesian product of complete graphs is the number of nontrivial factors.
\end{proof}

\begin{corollary}\label{cor:disjoint}
Let $\Delta$ be a pure $(d-1)$-dimensional simplicial complex with connected dual graph. If the minimal nonfaces $M_1,\dots,M_s$ of $\Delta$ are pairwise disjoint, then 
\[\Delta=\partial\langle M_1\rangle*\dots*\partial\langle M_s\rangle*\big\langle V\setminus\textstyle\bigcup M_i\big\rangle\] is pure and normal, $d=n-s$, and $\diam(\Delta)=s=n-d$. If moreover $\Delta$ is a sphere, then $V=\bigcup M_i$; if $\Delta$ is the boundary of a simplicial polytope $P$, the polar dual of $P$ is combinatorially a product of $s$ simplices, of dimensions $|M_1|-1,\dots,|M_s|-1$.
\end{corollary}

\begin{proof}
Apply Lemma~\ref{lem:join} with $W=V$ and $Z_i=M_i$ (each of size $\ge2$). Links are again of this form, on $W\setminus\sigma$ with the sets $M_i\setminus\sigma$, so they have connected dual graphs. If $V\ne\bigcup M_i$, $\Delta$ is a cone, hence not a sphere. Joins of simplex boundaries are dual to products of simplices.
\end{proof}

So for pairwise disjoint minimal nonfaces the Hirsch bound holds with equality. Main Theorem~\ref{thm:C} allows overlaps, confined to a fixed set $X$. 

\begin{lemma}\label{lem:fiber}
Under the hypotheses of Main Theorem~\ref{thm:C}, let $Y\subseteq X$ be a face of $\Delta$ and
\begin{gather*}
\Delta_Y:=\{\tau\subseteq V\setminus X:\ \tau\cup Y\in \Delta\},\\ \M_Y:=\{M \text{ minimal nonface}:\ M\cap X\subseteq Y\}.
\end{gather*}
Then the sets $M\setminus X$, $M\in\M_Y$, are nonempty and pairwise disjoint, and \[\Delta_Y=\{\tau\subseteq V\setminus X:\ M\setminus X\not\subseteq\tau\text{ for all }M\in\M_Y\}.\] Consequently $\Delta_Y$ is as in Lemma~\ref{lem:join} with $W=V\setminus X$ and $t=|\M_Y|$.
\end{lemma}

\begin{proof}
For $\tau\subseteq V\setminus X$, $\tau\cup Y\notin \Delta$ if and only if some minimal nonface $M$ satisfies $M\subseteq\tau\cup Y$, i.e.\ $M\cap X\subseteq Y$ and $M\setminus X\subseteq\tau$. If $M\in\M_Y$ had $M\setminus X=\emptyset$, then $M\subseteq Y$, contradicting $Y\in \Delta$. For distinct $M,M'$, $(M\setminus X)\cap(M'\setminus X)=(M\cap M')\setminus X=\emptyset$.
\end{proof}

\begin{proof}[\bf Proof of Main Theorem~\ref{thm:C}]
Write $V=[n]$. For a facet $F$ of $\Delta$, call $F\cap X$ its trace. We first show that any two facets with the same trace can be joined by a path of length at most $n-d$, all of whose facets have that same trace.
Fix a set $Y\subseteq X$ that occurs as a trace, and choose a facet $F$ with $F\cap X=Y$. Since $Y\subseteq F$, the set $Y$ is a face of $\Delta$, so Lemma~\ref{lem:fiber} applies. The set $F\setminus X$ is a facet of $\Delta_Y$: it belongs to $\Delta_Y$, and any strictly larger face $\tau$ of $\Delta_Y$ would give a face $\tau\cup Y$ of $\Delta$ strictly containing $F$.
By Lemmas~\ref{lem:join} and~\ref{lem:fiber}, every facet of $\Delta_Y$ has $\lvert V\setminus X\rvert-\lvert\mathcal M_Y\rvert$ vertices.
Applying this to $F\setminus X$, we obtain $d-\lvert Y\rvert = \lvert V\setminus X\rvert-\lvert\mathcal M_Y\rvert$.
Consequently, $\lvert\mathcal M_Y\rvert = n-d-\lvert X\setminus Y\rvert\le n-d$.
The same calculation shows that every facet $\tau$ of $\Delta_Y$ has $d-\lvert Y\rvert$ vertices. By the definition of $\Delta_Y$, the set $\tau\cup Y$ is a face of $\Delta$.
It has exactly $d$ vertices, so it is a facet of $\Delta$; its trace is $Y$. Thus $\tau\longmapsto\tau\cup Y$ is a bijection from the facets of $\Delta_Y$ to the facets of $\Delta$ with trace $Y$. Moreover, $\lvert(\tau\cup Y)\cap(\tau'\cup Y)\rvert
=
\lvert\tau\cap\tau'\rvert+\lvert Y\rvert$,
so this bijection preserves adjacency in both directions. The facets with trace $Y$ therefore induce a subgraph of $G(\Delta)$ isomorphic to $G(\Delta_Y)$. By Lemma~\ref{lem:join}, any two of them can be joined by a path of length at most $\lvert\mathcal M_Y\rvert\le n-d$ that stays entirely within this trace class.

We now connect two arbitrary facets $A$ and $B$.
Since $G(\Delta)$ is connected, choose a path $P_0=A,\ P_1,\ \ldots,\ P_\ell=B$.
Let $j_0$ be the last index for which $P_{j_0}$ has the same trace as $A$. Replace the portion $P_0,\ldots,P_{j_0}$ by a path of length at most $n-d$ whose facets all have this trace. Such a path exists by the argument above. If $j_0=\ell$, we have reached $B$.
Otherwise, keep the edge from $P_{j_0}$ to $P_{j_0+1}$ and repeat the construction on the remaining portion $P_{j_0+1},\ldots,P_\ell$ of the original path.
At each stage, we move to the last occurrence of the current trace in the portion of the original path under consideration.
No facet of the remaining portion has that trace. Each later replacement path stays within the trace class of its starting facet, so it cannot return to a class already used. The construction ends at $B$ and produces a path in which each trace class is visited in a single consecutive segment.
Suppose that this new path visits $\rho$ trace classes. Then $\rho\le\kappa$. Within each class the path has length at most $n-d$, and there are exactly $\rho-1$ edges between successive classes. Its total length is therefore at most
\[
\rho(n-d)+(\rho-1)
=
\rho(n-d+1)-1
\le
\kappa(n-d+1)-1.
\]
Since $A$ and $B$ were arbitrary, this proves the first inequality. Finally, every trace is a subset of $X$, so $\kappa\le 2^{\lvert X\rvert}=2^g$, which gives the second inequality.
\end{proof}

\begin{remark}
The proof gives slightly more: $\diam(\Delta)\le\sum_{Y}(n-d-|X\setminus Y|+1)-1$, the sum over the traces $Y$ met by the path. Since every trace is a face of the induced subcomplex $\Delta[X]$, $\kappa$ is at most the number of faces of $\Delta[X]$.
\end{remark}

\begin{remark} The conclusion of Main Theorem~\ref{thm:C} cannot be strengthened to the Hirsch bound  $\diam (\Delta) \le n-d$. To see this, let $\Delta$ be the pure $5$-dimensional complex on $\{0,\dots,7\}$ whose minimal nonfaces are $013,034,123,24$, so that $d=6$, $n=8$, $X=\{0,\dots,4\}$, $g=5$, and $5,6,7$ are cone points. Thus $\Delta$ has the four facets $012567$, $014567$, $023567$, $134567$. Note that $\Delta$ is not normal: the facets containing $3$ are $023567$ and $134567$, which are not adjacent. However, the dual graph of $\Delta$ is connected, and one has \[\diam(\Delta)=3>2 = n-d.\]  
\end{remark}

\subsection{Monomial ideals that are not radical}
Under the Stanley--Reisner correspondence, $\Delta$ corresponds to the radical monomial ideal $I=I_\Delta$. In this language:
\begin{enumerate}[label=\textup{(\roman*)}]
\item $n-d=\operatorname{ht}I$, and $\G(\Delta)$ is the graph on the minimal primes of $I$ (see \cite{BV});
\item $\Delta$ is flag if and only if $I$ is generated by quadrics, and $\Chat$ corresponds to the ideal $(I_2)$ generated by the quadratic part of $I$;
\item colon ideals $(I:m)$, where $m\notin I$ is a squarefree monomial, correspond to stars; the link is obtained by also removing the variables of $m$.
\end{enumerate}

Now let $I\subseteq R=k[x_1,\ldots,x_s]$ be a non-radical monomial ideal, and let $\operatorname{Pol}(I)\subseteq R^{\mathrm{pol}}$ be its polarization, which basically replaces each power of a variable by a product of distinct, additional variables. 
Polarization preserves height, degrees of the minimal generators, and graded Betti numbers. If $R/I$ is Cohen--Macaulay, then
$R^{\mathrm{pol}}/\operatorname{Pol}(I)$ is Cohen--Macaulay as well, so its Stanley--Reisner complex is pure and normal. This provides the hypotheses needed for our diameter estimates. We define the polarized dual graph of $I$ as
\[
G_{\mathrm{pol}}(I):=G(\operatorname{Pol}(I)).
\]
This should not be confused with the graph $G(I)=G(\sqrt I)$ on the minimal primes of
$I$ from \cite{BV}.
Before applying our bounds, linear minimal generators may be eliminated together with their variables. This leaves $G_{\mathrm{pol}}(I)$ unchanged and can only decrease the height, so the bounds remain valid when expressed using the original height $\operatorname{ht}I$.
Two consequences are immediate:

\begin{itemize}
\item If $I$ is generated by quadrics, then $\diam\G_{\mathrm{pol}}(I)\le\operatorname{ht}I$ by  \cite{AB14}. 
\item Write the minimal generators as $m_j=\prod_i x_i^{a_{ij}}$, with $j=1,\ldots,t$. Their polarizations are $P(m_j)=\prod_i\prod_{\ell=1}^{a_{ij}}x_{i,\ell}$. For each $i$, let $e_i$ be the second largest of the exponents $a_{i1},\ldots,a_{it}$, with multiplicities
retained. If $t=1$, set $e_i=0$.
The variable $x_{i,\ell}$ occurs in at least two polarized minimal generators precisely when $\ell\le e_i$. Consequently all intersections of distinct generator supports are contained in $X=\{x_{i,\ell}:1\le\ell\le e_i\}$, where $|X|=\sum_i e_i=:g$.
Eliminating linear minimal generators does not change this overlap set, since their variables occur in no other minimal generator. After that elimination, the remaining supports are exactly the minimal nonfaces of the polarized complex. Hence $X$ satisfies the hypothesis of Main Theorem \ref{thm:C}, which gives $\diam G_{\mathrm{pol}}(I) \leq 2^g(\operatorname{ht}I+1)-1.$
When $g=0$, the minimal generators are pairwise coprime. If they all have degree at least two, Corollary~\ref{cor:disjoint} gives
$\operatorname{diam}G_{\mathrm{pol}}(I)
=\operatorname{ht}I$,
so the bound is attained.

\end{itemize}

\Needspace{6\baselineskip}\section{Distance from normality: proof of Main Theorem~\ref{thm:D}}\label{sec:normal}

Throughout this section $\Delta$ is pure of dimension $d-1$ with connected dual graph, $1\le q\le d-1$, $r=d-q$, and
\[
a=\binom dq,\qquad b=\binom{d-1}{q-1},\qquad T=f_{q-1}(\Delta),\qquad \frac ab=\frac dq>1 .
\]
For a $q$-face $\sigma$, the link $K_\sigma=\lk_\Delta(\sigma)$ is pure of rank $r$; let $m_\sigma=f_0(K_\sigma)$. The facets containing $\sigma$ form the \emph{star} of $\sigma$, whose dual graph is that of $K_\sigma$. The idea is to cut the distance layers of $\Delta$ into blocks, each controlled by the full star of one $q$-face.

\begin{lemma}[Stars of $q$-faces]\label{lem:stars}
Assume that $K_\sigma$ has connected dual graph for every $q$-face $\sigma$, and let $C\ge1$.
\begin{enumerate}[label=\textup{(\alph*)}]
\item If $\diam K_\sigma\le C(m_\sigma-r)$ for every $q$-face $\sigma$, then $\diam(\Delta)\le\frac{3C}{b}(T-a)$.
\item If $\diam K_\sigma+1\le C(m_\sigma-r+1)$ for every $q$-face $\sigma$, then $\diam(\Delta)\le\frac{3C}{b}(T-a)+C-1<\frac{3C}{b}\,T$.
\end{enumerate}
\end{lemma}

\begin{proof}
Let $F_0,F_D$ be the endpoints of a diameter, $D=\diam(\Delta)$, and let $L_0,\dots,L_D$ be the distance layers from $F_0$.
For a $q$-face $R$, the facets containing $R$ induce a connected subgraph.
Since the distance from $F_0$ changes by at most one along an edge, the indices of layers containing a facet that
contains $R$ form an interval.

We first divide the layers into consecutive blocks.
Put $s_1=0$. Given $s_i$, let $\ell_i$ be the largest $t$ such that some facet of $L_{s_i}$ and some facet of $L_t$ share a $q$-face. Choose such a face $S_i$, and put $s_{i+1}=\ell_i+1$. Continuing until the last
layer gives blocks $B_i=\{L_{s_i},\dots,L_{\ell_i}\}$, $i=1,\dots,k$.
By the interval property and the choice of $\ell_{i+1}$, every $q$-face lies in at most two consecutive blocks; indeed, if $R$ occurred in blocks $i$ and $j$ with $j\ge i+2$, the interval property would force it to
occur in $L_{s_{i+1}}$ and in a layer after
$L_{\ell_{i+1}}$, contradicting the choice of $\ell_{i+1}$.
The selected face $S_i$ lies in no layer after $\ell_i$, so every facet containing $S_i$ lies in block $i-1$ or $i$.

The full star of $S_i$ controls the length of $B_i$ through its first and last layers. Choose facets $A_i\in L_{s_i}$ and $A_i'\in L_{\ell_i}$ containing $S_i$. Their distances from $F_0$ differ by $\ell_i-s_i$.
They can also be joined inside the star of $S_i$, whose dual graph is that of $K_{S_i}$. Therefore
\[
\ell_i-s_i \leq \dist_{G(\Delta)}(A_i,A_i')
\leq \diam K_{S_i},
\]
and hence $|B_i|=\ell_i-s_i+1\leq \diam K_{S_i}+1$. This argument requires only the two endpoint facets to contain $S_i$.

Let $t_i$ be the number of $q$-faces of the star of $S_i$, and put $m_i=m_{S_i}$. To bound $\sum_i t_i$, fix a $q$-face $R$.
If $j$ is the first block in which $R$ occurs, then it occurs only in block $j$ and possibly block $j+1$. If $R$ belongs to the star of $S_i$, some facet contains both $R$ and $S_i$.
Its block must lie in both $\{j,j+1\}$ and $\{i-1,i\}$, which forces $i\in\{j,j+1,j+2\}$. Thus each $q$-face is counted at most three times, and $\sum_i t_i\le3T$.
Assume first that $k\ge2$. Each of the $a$ $q$-faces of $F_0$ occurs only in the first block, by the choice of $\ell_1$, so it belongs to at most the selected stars of $S_1,S_2$.
Each of the $a$ $q$-faces of $F_D$ occurs only in the last two blocks, so it belongs to at most the selected stars of $S_{k-1},S_k$. These two collections are
disjoint: indeed, a common $q$-face would force $\ell_1=D$ and hence $k=1$. Therefore $2a$ faces are counted at most twice, giving $\sum_i t_i\le3T-2a$.
For $k=1$ the same inequality follows from
$t_1\le T$ and $T\ge a$.

We next relate $t_i$ to $m_i$. The star of $S_i$ is pure with connected dual graph and $q+m_i$ vertices. Order its facets along a rooted spanning tree of its dual graph, with each facet after its parent. Each subsequent facet introduces at most one new vertex. The first facet has $a$ $q$-faces, and each of the $m_i-r$ further vertices is introduced by a facet containing $b$ new $q$-faces. Hence $t_i\ge a+b(m_i-r)$. 
Put $\varepsilon=1$ in case~(a) and $\varepsilon=C$ in case~(b). The hypothesis on $K_{S_i}$ now gives
\[
|B_i|
\leq\diam(K_{S_i})+1
\leq\varepsilon+C(m_i-r)
\leq\varepsilon+\frac{C}{b}(t_i-a).
\]
Summing over the blocks, we obtain
\begin{align*}
D+1 &=\sum_i|B_i| 
\le k\varepsilon+\frac{C}{b} \left(\sum_i t_i-ka\right) \le k\varepsilon+\frac{C}{b} \bigl(3T-(k+2)a\bigr)\\
&=\frac{3C}{b}(T-a)+\varepsilon -(k-1)\left(\frac{Ca}{b}-\varepsilon\right) \le \frac{3C}{b}(T-a)+\varepsilon,
\end{align*}
because $a/b=d/q>1$ and $C\geq\varepsilon$.
Taking $\varepsilon=1$ proves~(a), and taking
$\varepsilon=C$ proves the first inequality in~(b). Finally, $\frac{3C}{b}(T-a)+C-1<\frac{3C}{b}T$, since $a/b>1$ and $C\ge1$.
\end{proof}

\begin{proof}[\bf Proof of Main Theorem~\ref{thm:D}]
Note that $\mu_q=(T-a)/b$. Let $q\ge\nu(\Delta)$. Fix a face $\sigma$ with $|\sigma|=q$, and put
$K=\lk_\Delta(\sigma)$ and $m=f_0(K)$. 
For every face $\tau$ of $K=K_\sigma$ (including $\emptyset$) the facets of $K$ containing $\tau$ correspond to the facets of $\Delta$ containing $\sigma\cup\tau$, which induce a connected subgraph since $|\sigma\cup\tau|\ge q$. So $K$ is normal. 
To bound the flag deficiency of $K$, recall that $\lk_K(\tau) = \lk_\Delta(\sigma\cup\tau)$. 
If $|\tau|\ge \max\{\fd(\Delta)-q,0\}$, then $|\sigma\cup\tau| = q+|\tau| \ge \fd(\Delta)$. 
Hence $\lk_\Delta(\sigma\cup\tau)$ is flag.
It follows that $\fd(K) \le \max\{\fd(\Delta)-q,0\}$.

If $r=1$, $K$ is a set of $m$ points, and $\diam K\le m-1$. If $r\ge2$, Main Theorem~\ref{thm:B} gives $\diam K\le2^{\fd^*(K)}(m-r)\le 2^h(m-r)$, since $\fd^*(K)\le\min\{\fd(K),r-2\}$. Lemma~\ref{lem:stars}(a) with $C=2^h$ gives $\diam(\Delta)\le3\cdot2^h\mu_q$. For the second term, Corollary~\ref{cor:normal} applied to $K$, which has $m\le n-q$ vertices, gives $\diam K+1\le(m-r+1)B_{r,\,m-r+1}\le B\,(m-r+1)$ with $B=B_{r,\,n-d+1}$, because $m-r+1\le n-d+1$. Lemma~\ref{lem:stars}(b) with $C=B$ gives $\diam(\Delta)\le 3B\mu_q+B-1<3B(\mu_q+1)$. Finally, the spanning-tree count in the proof of Lemma~\ref{lem:stars}, applied to $\Delta$, gives $T\ge a+b(n-d)$, that is, $\mu_q\ge n-d$; and $T\le\binom nq$ gives the upper bound for $\mu_q$.
\end{proof}

For flag complexes the exponent $h$ vanishes, and the number of $q$-faces can be bounded in terms of $n$ and $d$ alone:

\begin{corollary}\label{cor:flagq}
Let $\Delta$ be a flag pure $(d-1)$-dimensional complex on $n$ vertices with connected dual graph, and let $q\ge\nu(\Delta)$ with $2\le q\le d-1$. Then
\[
\diam(\Delta)\ \le\ \frac{3d}{q}\Big(\Big(\frac nd\Big)^{q}-1\Big).
\]
\end{corollary}

\begin{proof}
Since $\Delta$ is flag, $\fd(\Delta)=0$, and
therefore $h=0$ in Main Theorem \ref{thm:D}. Hence
\[
\diam(\Delta) \leq 3\mu_q(\Delta)
= \frac{3}{\binom{d-1}{q-1}}
\left(f_{q-1}(\Delta)-\binom dq\right).
\]
So it remains to bound $f_{q-1}(\Delta)$, the number of faces with $q$ vertices.

Let $\Gamma$ be the $1$-skeleton of $\Delta$. Because $\Delta$ is flag, its faces with $q$ vertices are exactly the
$q$-cliques of $\Gamma$. Moreover, $\Gamma$ contains no $(d+1)$-clique (indeed, such a clique would be a face of $\Delta$, contrary to $\dim(\Delta)=d-1$).

We use the following numerical consequence of Zykov's theorem \cite{Zyk49}, stated as equation (1) in \cite{AS24}: an $n$-vertex graph containing no $(d+1)$-clique has at most $\binom dq\left(\frac nd\right)^q$ cliques with $q$ vertices.
This bound is valid whether or not $d$ divides $n$, and applying it to $\Gamma$ we get that $f_{q-1}(\Delta) \leq \binom dq\left(\frac nd\right)^q$.

Substituting into the diameter estimate we get
\[
\diam(\Delta) \leq \frac{3}{\binom{d-1}{q-1}} \left[\binom dq\left(\frac nd\right)^q-\binom dq \right]
= 3\,\frac{\binom dq}{\binom{d-1}{q-1}} \left(\left(\frac nd\right)^q-1\right)
= \frac{3d}{q} \left(\left(\frac nd\right)^q-1\right),
\]
where the last equality uses that $\binom{d}{q} /\binom{d-1}{q-1}=d/q$.
\end{proof}

\begin{remark}
For $q=1$, Main Theorem~\ref{thm:B} and Corollary~\ref{cor:normal} are better than Main Theorem~\ref{thm:D} (which gives $3(n-d)$ in the flag case, worse than the Adiprasito--Benedetti bound), and for $q=d-1$ Lemma~\ref{lem:ridges} is better. Main Theorem~\ref{thm:D} and Corollary~\ref{cor:flagq} matter in between, for $2\le q\le d-2$.
\end{remark}

\Needspace{6\baselineskip}\section{Distance from normality: a conditional Main Theorem~\ref{thm:E}}\label{sec:conj}

In Section~\ref{sec:normal}, flagness enters only through Main Theorem~\ref{thm:B}. In this section we state a flag version of the base of the Kalai--Kleitman recursion, and show what it would imply. Recall that a $1$-layered family of $r$-sets is flag if no facet lies in an $(r+1)$-clique of the graph whose edges are the pairs of vertices contained in a common facet; the facets of a flag complex have this property, and it is inherited by ranges of layers and by links of faces. 
For a layered family $\mathcal F$ and a face $S$ of $\mathcal F$, the \emph{link} $\lk_{\mathcal F}(S)$ is the family of the sets $F\setminus S$ with $S\subseteq F\in\mathcal F$, whose layers are the nonempty sets $L_t^S=\{F\setminus S: F\in L_t,\ S\subseteq F\}$. If $\mathcal F$ is $1$-layered, so is $\lk_{\mathcal F}(S)$. The \emph{flag deficiency} $\fd(\mathcal F)$ is the least $p$ such that $\lk_{\mathcal F}(S)$ is flag for every face $S$ with $|S|=p$; since links of flag families are flag, it is then flag for every face with at least $p$ vertices. If $\mathcal F$ consists of facets of a complex $\Delta$, then $\fd(\mathcal F)\le\fd(\Delta)$.

\begin{conjecture}\label{conj:flag}
Let $\mathcal F$ be a flag $1$-layered family of $r$-subsets of an $m$-set. Then $\mathcal F$ has at most $m-r+1$ layers.
\end{conjecture}

Adiprasito and Benedetti's theorem \cite{AB14} is the case in which $\mathcal F$ is the distance layering of a normal flag complex, so Conjecture~\ref{conj:flag} asks for the same bound for every layering. Flagness cannot be dropped: already for $r=2$ there are $1$-layered families of edges of graphs with triangles that have more than $m-1$ layers. The conjecture is trivial for $r=1$, and open for $r=2$, that is, for triangle-free graphs, where the Larman--Barnette bound only gives $2m-3$ layers (Theorem~\ref{thm:LBKK}). We tested it by computer on tens of thousands of random flag layered families with $r=2,3,4$, including distance layerings of random flag complexes, and found no counterexample.

\begin{lemma}\label{lem:flagclaim}
Assume Conjecture~\ref{conj:flag}. For $r\ge1$, every $1$-layered family of $r$-sets on $m\ge r$ vertices with flag deficiency at most $p$ has a number of layers not larger than 
\[(m-r+1)\binom{\lfloor\log_2(m-r+1)\rfloor+p}{p}.\]
\end{lemma}

\begin{proof}
We use the same decomposition as in Section~\ref{sec:LBKK}, but now we stop taking links when the family becomes flag.
Passing to a vertex link lowers the flag-deficiency bound $p$ by one and preserves the quantity $m-r+1$ used in the estimate. Passing to either outer range retains the bound $p$ and reduces that quantity to at most half its value. Conjecture~\ref{conj:flag} supplies the bound $m-r+1$ when $p=0$. This explains why the binomial factor involves $p$ in place of $r-1$.
We use induction on $p$ and on $m$.
Put $\mu=m-r+1$ and $\ell=\lfloor\log_2\mu\rfloor$. If $r=1$ there are at most $m=\mu$ layers, and if $m=r$ there is one. If $r\ge2$ and $p=0$ the family is flag, and there are at most $\mu$ layers by Conjecture~\ref{conj:flag}. Let $r\ge2$, $p\ge1$ and $m>r$, so that $\ell\ge1$. As in the proof of~\eqref{eq:rec}, the layers split into an initial range using at most $a=\lfloor(m-1)/2\rfloor$ vertices, a final range using at most $b=\lceil(m-1)/2\rceil$ vertices, and a middle range that is a range of $\lk_{\mathcal F}(v)$ for some vertex $v$, a $1$-layered family of $(r-1)$-sets on at most $m-1$ vertices. Ranges and links do not increase the flag deficiency, except that $\fd(\lk_{\mathcal F}(v))\le p-1$. Since $(m-1)-(r-1)+1=\mu$, the middle range has at most $\mu\binom{\ell+p-1}{p-1}$ layers by induction. For $x\in\{a,b\}$ with $x\ge r$ we have $x-r+1\le\mu/2$, so $\lfloor\log_2(x-r+1)\rfloor\le\ell-1$, and the corresponding range has at most $\frac\mu2\binom{\ell+p-1}{p}$ layers; if $x<r$ it is empty. Adding up and using Pascal's rule gives $\mu\binom{\ell+p}{p}$.
\end{proof}

\begin{proof}[\bf Proof of Main Theorem~\ref{thm:E}]
(a) The distance layering from a facet is a $1$-layered family of $d$-sets on $n$ vertices with $e+1$ layers, and its flag deficiency is at most $p$. By Lemma~\ref{lem:flagclaim} with $r=d$ and $m=n$, we have 
\[
\diam(\Delta)\le e<(n-d+1)\binom{\lfloor\log_2(n-d+1)\rfloor+p}{p}.
\]

(b) Put $r=d-q$. For a $q$-face $\sigma$, the link $K=\lk_\Delta(\sigma)$ is normal of rank $r$, as in the proof of Main Theorem~\ref{thm:D}, with $m\le n-q$ vertices and $\fd(K)\le p'$: a face of $K$ with $p'$ vertices, together with $\sigma$, has at least $p$ vertices. The distance layering of $K$ is a $1$-layered family of $r$-sets on $m$ vertices with flag deficiency at most $p'$, so by Lemma~\ref{lem:flagclaim},
\[
\diam K+1\le(m-r+1)\binom{\lfloor\log_2(m-r+1)\rfloor+p'}{p'}\le C\,(m-r+1),
\]
where $C=\binom{\lfloor\log_2(n-d+1)\rfloor+p'}{p'}$. Lemma~\ref{lem:stars}(b) gives $\diam(\Delta)<3C(\mu_q(\Delta)+1)$.
\end{proof}

\Needspace{6\baselineskip}\section{Examples}\label{sec:examples}
We give three small non-flag normal examples, one for each of Main Theorems~\ref{thm:A}, \ref{thm:B} and~\ref{thm:C}, to compare our bounds with those of Larman--Barnette, Novik and Kalai--Kleitman. The examples below were checked (and in the case of Section~\ref{sec:MeFamily}, generated) using Claude Sonnet 5.5. A fourth example, in Section~\ref{sec:exnormal}, is not normal and illustrates Main Theorem~\ref{thm:D}.

\subsection{Cross-polytopes minus a facet}\label{sec:exA}
Let $P_d$ be the boundary of the $d$-dimensional cross-polytope. It is a flag sphere on $n=2d$ vertices whose minimal nonfaces are the $d$ antipodal pairs. Let $\Delta_{2d}$ be obtained from $P_d$ by removing a single facet $F$ while retaining all its proper faces. 
For every proper face $\sigma$ of $\Delta_{2d}$, the facets containing $\sigma$ induce a cube, or a cube with one vertex deleted, according as $\sigma$ is not contained in $F$ or is contained in $F$.
These graphs are connected, and facet stars are singletons. Hence $\Delta_{2d}$ is normal. Its minimal nonfaces are the $d$ antipodal pairs and the set $F$, which has size $d$. Hence $\Delta_{2d}$ satisfies $(\mathrm G_1)$ for $d\ge 3$, and $\Chat=P_d$ is pure and normal of dimension $d-1$. The link of a face with $d-2$ vertices is a $4$-cycle in $P_d$, and a path with three edges when the face lies in $F$. In both cases the dual graph has diameter $2$, so $\ldiam_1(\Delta_{2d})=2$ and Main Theorem~\ref{thm:A} gives
\[
\diam(\Delta_{2d})\le (n-d)\cdot\ldiam_1=2d.
\]
The dual graph of $P_d$ is the $d$-cube $Q_d$, so $\G(\Delta_{2d})$ is $Q_d$ minus a vertex $v$. Identify the vertices of $Q_d$ with the subsets of $[d]$ and $v$ with $\emptyset$. Between two nonempty sets $A,B$, first add the elements of $B\setminus A$ and then remove those of $A\setminus B$; this is a Hamming geodesic that avoids $\emptyset$. So $\G(\Delta_{2d})$ is isometrically embedded in $Q_d$, and complementary nonempty sets show that $\diam(\Delta_{2d})=d$ for all $d\ge3$. On the other hand, $\fd^{*}(\Delta_{2d})=d-2$ and $\beta(\Delta_{2d})=d$: the face $F$ is a critical clique with $d$ vertices that is not a face. Hence Main Theorem~\ref{thm:B} only gives $2^{d-2}(n-d)=2^{d-2}d$, which is the Larman--Barnette bound, and Novik's bound is $2^{d-2}n$.

\subsection{A locally flag family}\label{sec:MeFamily}
The inequality of Proposition~\ref{prop:banner} can be as far from equality as possible. We index the examples by their number $n$ of vertices and by the number $d$ of vertices of their facets.

For a sequence $s=(s_1,\ldots,s_N)$ of distinct vertices and $1\le m\le N$, let
\[
W(s;m)=
\left\langle
\{s_j,\ldots,s_{j+m-1}\}:1\le j\le N-m+1
\right\rangle.
\]
Thus its facets are the sets of $m$ consecutive terms of $s$. This complex is flag: the positions occupied by a clique differ by at most $m-1$, so the clique lies in a facet. It is also normal, since the facets containing a given face have consecutive starting indices. Its dual graph is a path when $m\ge2$ and is complete when $m=1$.

\begin{definition}\label{def:Me}
Fix $d\ge3$ and $n\ge d^2-d+2$, and set $L:=n-(d^2-d+2)\ge0$. The vertex set of $M_{n,d}$ consists of $X=\{x_1,\dots,x_{d}\}$ (indices mod $d$), the vertices $w$ and $r$, the sets $U_i=\{u_{i,1},\dots,u_{i,d-2}\}$ for $i\in[d]$, and $L$ further vertices $v_1,\dots,v_L$; so there are $d+2+d(d-2)+L=n$ vertices. Let
\begin{gather*}
S_i=(x_{i+1},\dots,x_{i+d-1},u_{i,1},\dots,u_{i,d-2},w),\\ K_i=x_i * W(S_i;d-1),
\end{gather*}
let $R=W(U_1U_2\cdots U_{d}v_1\cdots v_L;\,d-2)$ (on the concatenated sequence), and
\[
M_{n,d}:=K_1\cup\dots\cup K_{d}\cup\big(\{w,r\}*R\big).
\]
\end{definition}

When $d=3$, every facet of $R$ consists of one vertex, so its dual graph is complete.
For $d=3$ and $L=0$, with $x_i=i$, $w=4$, $r=5$, $u_{1,1}=6$, $u_{2,1}=7$, $u_{3,1}=8$, the complex $M_{8,3}$ is the complex on $8$ vertices with facets
\[123,\,136,\,146,\,127,\,247,\,238,\,348,\,456,\,457,\,458.\]

For $d\geq 4$, the construction has a simple description in the dual graph. Each $K_i$ contributes a path of $d-1$ edges from the common facet $X$ to $A_i=\{x_i,w\} \cup U_i$. The complex $\{w,r\}*R$ contributes a further path.
Along this path, the facets $B_i=\{w,r\}\cup U_i$ occur at intervals of $d-2$ edges, and $A_i$ is adjacent to $B_i$. Thus the path from $X$ to $B_i$ has length $d$. Appending $v_1,\ldots,v_L$ adds a terminal path of $L$ edges after $B_d$.

\begin{remark}
\label{rem:banner num Mnd}
The set $X\cup\{w\}$ is a critical clique of size $d+1$: removing $w$ leaves the facet $X$, and $w$ is adjacent to every $x_i$. Since $M_{n,d}$ has dimension $d-1$, this clique is not a face. Consequently the banner number of $M_{n,d}$ is $\beta(M_{n,d})=d+1$.
\end{remark}

\begin{proposition}\label{prop:Me}
For every $d\ge3$ and $n\ge d^2-d+2$, $M_{n,d}$ is a pure normal $(d-1)$-dimensional complex on $n$ vertices with $\beta(M_{n,d})=d+1$ and $\fd(M_{n,d})=1$. Moreover:
\begin{enumerate}[label=\textup{(\alph*)}]
\item for $d\ge6$, $\diam(M_{d^2-d+2,\,d})=2\lceil 3(d-1)/2\rceil=3\sqrt n+O(1)$;
\item for $d\ge5$ and $n\ge d^2-d+2+\lfloor(d-2)/2\rfloor$, $\diam(M_{n,d})=n-d^2+3d-2+\lfloor(d-2)/2\rfloor$;
\item for $d=4$ and $n\ge14$, $\diam(M_{n,4})=n-7$.
\end{enumerate}
\end{proposition}

Thus Novik's bound for $M_{n,d}$ is $2^{d-1}n$, while Main Theorem~\ref{thm:B} gives $2(n-d)$, and for $d\ge6$ the true diameter is of order $\sqrt n$ when $n=d^2-d+2$ and of order $n$ (with slope $1$) for large $n$. For $d=4$ the complexes $M_{n,4}$ miss the Hirsch value $n-4$ by exactly $3$. The complex $M_{14,4}$ has $14$ vertices and $M_{22,5}$ has $22$; a computer check gives $\diam(M_{14,4})=7$ and $\diam(M_{22,5})=11$, with $\fd^{*}=1$ in both cases. For $d=3$ the diameter is $3$ if $n=8$ and $4$ for every $n\ge9$.

\begin{proof}[Proof of Proposition \ref{prop:Me}]
Put $M=M_{n,d}$, $B=\langle\{w,r\}\rangle*R$, and $A_i=\{x_i,w\}\cup U_i$.
All generating facets have $d$ vertices, and all $n$ prescribed vertices occur. Thus $M$ is pure of dimension $d-1$. Remark~\ref{rem:banner num Mnd} gives $\beta(M)=d+1$ and shows that $M$ is not flag.

We first check that every vertex link is flag. The complexes $K_i$ and $B$ are flag, as are their links. We use the fact that a union of two flag complexes whose intersection is a full simplex is flag. Indeed, no edge joins vertices exclusive to different pieces, so every clique lies in one piece.

Since $w$ is the last entry of $S_i$, the only set of $d-1$ consecutive terms containing $w$ is
$U_i\cup\{w\}$. Hence $A_i$ is the only facet of $K_i$ containing $w$. We have $\lk_M(r)=\langle w\rangle*R$ and $\lk_M(v_a) = \langle\{w,r\}\rangle*\lk_R(v_a)$, which are flag. Also, $\lk_M(w) = (\langle r\rangle*R) \cup \bigcup_{i=1}^d\langle\{x_i\}\cup U_i\rangle$.
Each added simplex meets the preceding union in $\langle U_i\rangle$, so this link is flag.
For $u\in U_i$, $\lk_M(u) = \lk_{K_i}(u)\cup\lk_B(u)$, and the two pieces intersect in the full simplex $\langle\{w\}\cup(U_i\setminus\{u\})\rangle$. Finally, $\lk_M(x_i) = \bigcup_{j=1}^d\lk_{K_j}(x_i)$.
Any two pieces meet in the full simplex $\langle X\setminus\{x_i\}\rangle$: the sets $U_j$ are disjoint, and $w$ belongs only to the piece coming from $K_i$. Successive applications of the same gluing observation prove flagness of these links. Since $M$ is not flag, $\fd(M)=1$.

For normality, we use another gluing observation. If two pure normal complexes of dimension $t$ intersect in a pure $(t-1)$-dimensional complex, then their union is normal. Indeed, a face in the intersection is contained in a common $(t-1)$-face.
Facets extending that face in the two complexes are adjacent, joining the two connected stars. A face outside the intersection occurs in only one of the two complexes.

Let $K_i'$ be the subcomplex generated by the facets of $K_i$ other than $A_i$, and put $B'=\bigcup_{i=1}^d K_i'$. Starting with $X$, add the remaining facets in their order along each sequence $S_i$. Each added facet introduces a new vertex of $U_i$ and meets the preceding union in exactly a ridge. Hence $B'$ is normal. The complex $L'=\lk_M(w)$ is also normal: start with $\langle r\rangle*R$ and attach each simplex $\langle\{x_i\}\cup U_i\rangle$ along the ridge $\langle U_i\rangle$. Now 
$
M=B'\cup(\langle w\rangle*L'),
$
and 
$
B'\cap(\langle w\rangle*L')
=
\bigcup_{i=1}^d\langle\{x_i\}\cup U_i\rangle$.
The intersection is pure of dimension $d-2$. Since cones over normal complexes are normal, the gluing observation proves that $M$ is normal.

We now compute the diameter for $d\ge4$. Let $z_1,\ldots,z_N$ be the concatenation $U_1U_2\cdots U_dv_1\cdots v_L$, and put $H_j=\{w,r,z_j,\ldots,z_{j+d-3}\}$, with $1\le j\le\ell+L$ and $\ell=(d-1)(d-2)+1$. These facets induce the path $H_1,\ldots,H_{\ell+L}$.
For each $i$, the facets of $K_i$ give a path of $d-1$ edges from $X$ to $A_i$, and $A_i$ is adjacent to
$
H_{q_i}=\{w,r\}\cup U_i$,
with $q_i=(i-1)(d-2)+1$.
Denote the resulting path of length $d$ from $X$ to $H_{q_i}$ by $P_i$.

The paths $P_i$ have no common internal vertices, and their internal vertices do not lie on the path of the $H_j$. There are no further edges: facets from distinct $K_i$, other than $X$, share at most $d-2$ vertices, and the only adjacency between $K_i$ and $B$ is $A_iH_{q_i}$.
This follows because every other facet of $K_i$ contains neither $w$ nor $r$, while $A_i$ shares $d-1$ vertices with a facet of $B$ precisely when that facet contains all of $U_i$.

Write $\mu(j)=\min_{1\le i\le d}|j-q_i|$, and $m=\left\lfloor\frac{d-2}{2}\right\rfloor$.
A shortest path between two facets $H_j,H_{j'}$
either stays on the path of the $H_j$ or passes
through $X$. Therefore
\[
\begin{aligned}
\operatorname{dist}(X,H_j)
&=d+\mu(j),\\
\operatorname{dist}(H_j,H_{j'})
&=
\min\{|j-j'|,\;2d+\mu(j)+\mu(j')\}.
\end{aligned}
\]

First take $L=0$, and write $G_0$ for the resulting dual graph. Vertices on the paths $P_i$ are at distance at most $d$ from $X$, and $\mu(j)\le m$ for $1\le j\le\ell$. Hence $\diam G_0\le2d+2m$. For $d\ge6$, the facets $H_{m+1}$ and $H_{\ell-m}$ both have $\mu$-value $m$. Moreover, $\ell-1=(d-1)(d-2)\ge4d-4\ge2d+4m$. Thus their distance is $2d+2m$, and
\[
\diam G_0 = 2d+2m = 2\left\lceil\frac{3(d-1)}2\right\rceil.
\]
Since $n=d^2-d+2$ when $L=0$, this is also
$3\sqrt n+O(1)$, proving~(a).

For $d\ge5$, every vertex of $G_0$ is at distance at most $2d+m$ from $H_\ell$, by going through $X$. This bound is attained at $H_{m+1}$, because
\[
\dist(H_\ell,H_{m+1}) =
\min\{\ell-1-m,\;2d+m\} = 2d+m.
\]
Here $\ell-1\ge2d+2m$: equality holds when $d=5$,
and for $d\ge6$ it follows from the inequality above.

For general $L$, the additional facets form the path $H_\ell,H_{\ell+1},\ldots,H_{\ell+L}$, which meets $G_0$ only at $H_\ell$. Consequently, $\diam(M_{n,d})= \max\{\diam G_0,\;L+2d+m\}$. 
The hypothesis in~(b) is exactly $L\ge m$. Since $\diam G_0\le2d+2m$, we obtain
\[
\diam(M_{n,d}) =L+2d+m =n-d^2+3d-2+
\left\lfloor\frac{d-2}{2}\right\rfloor,
\]
proving~(b).

For $d=4$, first take $L=0$.
Now $\ell=7$, $q_i=2i-1$, and each $P_i$ has length $4$. Two facets on the same $P_i$ are at distance at most $4$, while two facets among $H_1,\ldots,H_7$ are at distance at most $6$.
For facets $F\in P_i$ and $F'\in P_k$ with $i\ne k$, put $u=\dist(X,F)$ and $v=\dist(X,F')$.
Considering paths through $X$ and paths avoiding
$X$ gives $\dist(F,F') = \min\{u+v,\;8-u-v+|q_i-q_k|\}\le7$, since the two quantities have sum at most $14$.
Similarly, $\dist(F,H_j) = \min\{u+4+\mu(j),\;4-u+|q_i-j|\} \le7$.
For the last inequality, use $\mu(j)\le\min\{j-1,7-j\}$ and $|q_i-j|\le\max\{j-1,7-j\}$; again the two quantities have sum at most $14$.
Thus $\diam G_0\le7$.
The facet $A_1$ lies at distance $3$ from $X$ on $P_1$, so the preceding distance formula gives $\dist(A_1,H_7)=\min\{3+4,\;1+6\}=7$.
Therefore $\diam G_0=7$, and the largest distance from $H_7$ is also $7$. Appending the path of $L$ edges gives diameter $L+7=n-7$, proving~(c).

For completeness, when $d=3$, the facets of $B$ are pairwise adjacent, and there are three paths of length $3$ from $X$ to $B_i=\{w,r\}\cup U_i$.
When $L=0$, two facets on different such paths, at distances $u,v$ from $X$, can be joined by routes of lengths $u+v$ and $7-u-v$.
Their minimum is at most $3$, because these are integers. Facets on the same path are also at distance at most $3$, and $\dist(X,B_i)=3$.
Thus $\operatorname{diam}(M_{8,3})=3$.

If $L\ge1$, each additional facet $\{w,r,v_a\}$ is adjacent to all the $B_i$ and is at distance $4$ from $X$. From a facet at distance $u$ from $X$ on the path to $B_i$, it can be reached in at most $4-u$ steps. The additional facets are pairwise adjacent, so the diameter is $4$.
\end{proof}

\subsection{A complex with almost disjoint nonfaces}\label{sec:exC}
Let $d=4$. Let $\Gamma_7$ be the $(d-1)$-dimensional complex on the seven vertices $\{0,\dots,6\}$, whose minimal nonfaces are $01$, $04$, $23$ and $456$.  
In other words, 
\[
\Gamma_7 = 0256,\ 0356,\ 1245,\ 1246,\ 1256,\ 1345,\ 1346,\ 1356.
\]
The only overlaps among non-faces are $01\cap 04=\{0\}$ and $04\cap 456=\{4\}$, so $X=\{0,4\}$ and $g=2$. The traces $F\cap X$ are $\emptyset$, $\{0\}$ and $\{4\}$ (the set $\{0,4\}$ is a nonface), so $\kappa=3$. Main Theorem~\ref{thm:C} gives
\[
\diam(\Gamma_7)\le \kappa\,(n-d+1)-1=3\cdot 4-1=11.
\]
In fact $\Gamma_7=\partial\langle23\rangle*C$ with $C=\langle056,145,146,156\rangle$.
The last three facets of $C$ form a cone over the triangle graph $45,46,56$, and the facet $056$ is attached along the edge $56$.
The dual graph of $C$ consists of the triangle on $145,146,156$, with the vertex $056$ attached to $156$. The facets containing each vertex induce a connected subgraph: for vertices $5$ and $6$ these are the paths $056,156,145$ and $056,156,146$, respectively; the remaining vertex stars induce a triangle, an edge, or a single vertex. Stars of edges and facets are also connected, so $C$ is normal, and so is the join $\Gamma_7$.
The dual graph of $C$ is a triangle with a pendant vertex, of diameter $2$, and $\G(\Gamma_7)=K_2\,\square\,\G(C)$. Hence $\Gamma_7$ is normal with $\diam(\Gamma_7)=3=n-d$, so it satisfies the Hirsch bound. Here $\beta(\Gamma_7)=5=d+1$, since $\{1,2,4,5,6\}$ is a critical clique (removing $4$ leaves the facet $1256$) that is not a face. So Novik's bound is $2^{3}n=56$. Also $\fd^{*}(\Gamma_7)=2$, and Main Theorem~\ref{thm:B} gives $2^2(n-d)=12$. Main Theorem~\ref{thm:A} does not apply since the flag completion $\hat{\Gamma_7}$ is not $3$-dimensional.
Joining $\Gamma_7$ with disjoint boundaries of triangles keeps $\kappa=3$: each such factor adds $3$ to $n$, $2$ to $d$ and $1$ to the diameter, so the bound of Main Theorem~\ref{thm:C} grows by $3$ per factor, while the Larman--Barnette bound is multiplied by roughly $4$.

\subsection{Comparison of the previous examples}\label{sec:comparison}
Table~\ref{tab:examples} compares the three examples. The Larman--Barnette bound is $2^{d-2}(n-d)$ (the case $\fd^{*}=d-2$ of Main Theorem~\ref{thm:B}), Novik's is $2^{\beta-2}n$, and for Kalai--Kleitman we use the sharper form $(n-d+1)\binom{\lfloor\log_2(n-d+1)\rfloor+d-1}{d-1}-1$ of Corollary~\ref{cor:normal}, instead of the original $n^{\log_2 d+1}-1$ \cite[Thm.~3.12]{San13}, which is larger.

\begin{table}[h]
\setlength{\tabcolsep}{3.5pt}
\begin{tabular}{lrrrrrrr}
Example & $n$ & $d$ & $\diam$ & Our bound & Larman & Novik & Kalai--Kleitman\\
\hline
$\Delta_{12}$ (Sec.~\ref{sec:exA}) & $12$ & $6$ & $6$ & \phantom{22}$12$ (Thm A) & $96$ & $192$ & $146$\\
$M_{22,5}$ (Sec.~\ref{sec:MeFamily}) & $22$ & $5$ & $11$ & \phantom{22}$34$ (Thm B) & $136$ & $352$ & $1259$\\
$\Gamma_7$ (Sec.~\ref{sec:exC}) & $7$ & $4$ & $3$ & \phantom{22}$11$ (Thm C) & $12$ & $56$ & $39$\\ \medskip
\end{tabular}
\caption{The three examples.} 
\label{tab:examples}
\end{table}

In all three cases one of our bounds is smaller than all the others. For the families, the comparison improves with $d$: for $\Delta_{2d}$ our bound is $2d$ and the diameter is $d$, whereas $2^{d-2}d$ and the Kalai--Kleitman bound $(d+1)\binom{\lfloor\log_2(d+1)\rfloor+d-1}{d-1}-1$ grow much faster; for $M_{n,d}$ our bound is $2(n-d)$ and Novik's is $2^{d-1}n$. The Kalai--Kleitman expression is smaller than the Larman--Barnette bound only when $d$ is large compared with $\log n$.

Table~\ref{tab:examples} is in the regime where $d$ is small compared with $\log n$, which is unfavorable to Kalai--Kleitman. For the cross-polytope examples $\Delta_{2d}$ of Section~\ref{sec:exA}, which have $n=2d$, the Kalai--Kleitman bound is smaller than the Larman--Barnette bound for every $d\ge10$ that we checked ($d\le400$). Table~\ref{tab:larged} shows two larger values of $d$, and shows that in this regime too our theorem gives better bounds.

\begin{table}[h]
\setlength{\tabcolsep}{3.5pt}
\begin{tabular}{lrrrrrrr}
Example & $n$ & $d$ & $\diam$ & Thm A & Larman & Novik & Kalai--Kleitman\\
\hline
$\Delta_{80}$ & $80$ & $40$ & $40$ & $80$ & $1.1\cdot10^{13}$ & $2.2\cdot10^{13}$ & $4.5\cdot10^{7}$\\
$\Delta_{200}$ & $200$ & $100$ & $100$ & $200$ & $3.2\cdot10^{31}$ & $6.3\cdot10^{31}$ & $1.6\cdot10^{11}$\\ \medskip
\end{tabular}
\caption{Cross-polytopes minus a facet, in large dimension. ``Kalai--Kleitman'' refers to the sharper form of the Kalai--Kleitman theorem proved in Corollary~\ref{cor:normal}.}
\label{tab:larged}
\end{table}

\subsection{A complex that is not normal}\label{sec:exnormal}
Let $n=30$ and $d=6$. Call a sequence $v_1,\dots,v_N$ of elements of $\{0,\dots,n-1\}$ \emph{diverse} if the pairs $\{v_i,v_j\}$ with $1\le|i-j|\le5$ are pairwise distinct (in particular $v_i\ne v_j$ for such $i,j$). A diverse sequence has $5N-15$ such pairs, so $N\le(\binom n2+15)/5=90$. A randomized computer search found the following diverse sequence with $N=57$:
\[
\begin{array}{c}
0,\,24,\,26,\,11,\,4,\,13,\,1,\,20,\,6,\,28,\,29,\,19,\,16,\,8,\,14,\,9,\,13,\,27,\,17,\\
21,\,12,\,18,\,11,\,10,\,28,\,2,\,7,\,5,\,15,\,27,\,6,\,26,\,23,\,25,\,22,\,17,\,16,\,2,\\
1,\,0,\,3,\,21,\,8,\,23,\,7,\,24,\,20,\,12,\,14,\,25,\,5,\,3,\,4,\,18,\,9,\,22,\,15.
\end{array}
\]
Let $\Omega_{30}$ be the pure $5$-dimensional complex on $\{0,\dots,29\}$ whose facets are the $52$ sets $T_i=\{v_i,v_{i+1},\dots,v_{i+5}\}$, $1\le i\le52$.

Suppose that $i\ne j$ and that $T_i,T_j$ share two vertices $u,w$. By the defining property of the sequence, the two occurrences of the pair $\{u,w\}$ use the same two positions. Those positions therefore belong to both index intervals $[i,i+5]$ and $[j,j+5]$.
In particular, $T_i$ and $T_j$ cannot be equal: applying this argument to all their pairs would place six distinct positions in the intersection of the two distinct index intervals, which has at most five positions.
If $T_i,T_j$ share five vertices, the same argument applied to their pairs shows that the five corresponding positions belong to both index intervals. Since $i\ne j$, this forces $|i-j|=1$.
Conversely, consecutive facets share exactly five vertices. Thus $G(\Omega_{30})$ is the path
$T_1,\ldots,T_{52}$, and $\diam(\Omega_{30})=51$.
If the positions of an edge $\{u,w\}$ are $a<b$, the facets containing it are exactly the $T_i$ with $\max\{1,b-5\}\le i\le\min\{a,52\}$.
Their indices form an interval. For a face with at least two vertices, the facets containing it form the intersection of the intervals associated with its edges. This is again an interval, so $\Omega_{30}$ is $2$-normal.
On the other hand, the label $0$ occurs only at positions $1$ and $40$ (a direct check). Hence the facets containing $0$ are exactly $T_1$ and $T_{35},\ldots,T_{40}$. They induce two components in the dual path, so $\Omega_{30}$ is not normal. Therefore $\nu(\Omega_{30})=2$.
Since $d-1=5$, this is not the extreme case $q=d-1$.

Since $\nu(\Omega_{30})\ge2$, the Larman--Barnette, Novik and Kalai--Kleitman bounds do not apply, because $\Omega_{30}$ is not normal. The bound of Lemma~\ref{lem:ridges} yields
\[
\diam(\Omega_{30})\le\Big\lfloor\frac{\binom{30}{5}-6}{5}\Big\rfloor=28500 .
\]
Main Theorem~\ref{thm:D} with $q=2$ and $r=4$ gives instead a much better bound. A direct check shows that $\fd(\Omega_{30})=1$, so $h=\max\{0,\min\{\fd-q,\,r-2\}\}=0$. Moreover, the edges of $\Omega_{30}$ are exactly the $5N-15=270$ pairs counted above, so $f_1(\Omega_{30})=270$ and $\mu_2(\Omega_{30})=(270-15)/5=51$. Hence
\[
\diam(\Omega_{30})\le3\,\mu_2(\Omega_{30})=153,
\]
which is three times the real value $51$ of the diameter, but about $186$ times smaller than the previous bound. 

\section{Questions}\label{sec:questions}

\begin{question}\label{q:normalhat}
In Main Theorem~\ref{thm:A} with $k\ge2$, does normality of $\Chat$ follow from normality of $\Delta$ together with $(\mathrm G_k)$ and purity of $\Chat$? (For $k\le1$ it does, by Proposition~\ref{prop:normalhat}.)
\end{question}

\begin{question}
In Main Theorem~\ref{thm:C}, can $2^g$ be replaced by a polynomial in $g$ when $\Delta$ is normal? Note that every complex satisfies the hypothesis with $X=V$ and $g=n$, so a bound of the form $\mathrm{poly}(g)\,(n-d+1)$, or even $n-d+O(g)$, would give a polynomial, respectively linear, bound for all normal complexes. Is $\diam(\Delta)\le n-d$ for normal $\Delta$ when $g$ is small? (This cannot hold for all $g$, since Santos' counterexamples are normal. A random search over a few hundred normal complexes with small $g$ and $d$ found no complex with $\diam(\Delta)>n-d$.)
\end{question}

\begin{question}
Do normal complexes with $\fd(\Delta)=1$ satisfy the Hirsch bound? By Main Theorem~\ref{thm:B} and Proposition~\ref{prop:Me}, the maximum of $\diam(\Delta)$ over normal $3$-dimensional complexes with $\fd=1$ on $n$ vertices lies between $n-7$ (for $n\ge14$) and $2(n-4)$.
\end{question}

\begin{question}\label{q:subflag}
Let $\Delta$ be a normal subflag complex of dimension $d-1$ on $n$ vertices, that is, the $(d-1)$-skeleton of a flag complex. Is $\diam(\Delta)\le n-d$, or at least bounded by a polynomial in $n$?
\end{question}

\begin{question}[Direct interpolation]\label{q:direct}
Let $h\ge1$ and let $\Delta$ be a pure normal complex on $n$ vertices all of whose minimal nonfaces have at most $h+1$ vertices (so $h=1$ is the flag case). Is $\diam(\Delta)\le n^{h}$, or even $\diam(\Delta)\le c_h\,n$ with $c_h$ independent of the dimension? 
\end{question}

Main Theorem~\ref{thm:A} goes in the opposite direction: it bounds $\diam(\Delta)$ for complexes whose minimal nonfaces are either edges or \emph{large}.

\begin{question}[{\cite[Question D]{LMS17}}]
If every minimal nonface of a pure normal complex $\Delta$ has at most $k$ vertices, is $\diam(\Delta)\le 2^{k-2}n$? (We read ``size'' where \cite{LMS17} says ``dimension'', so that $k=2$ is the flag case, as for $k$-banner complexes.) Note that $\fd(\Delta)$ is not bounded in terms of the size of the minimal nonfaces (a cone over a non-flag complex has $\fd$ as large as one likes), so Main Theorem~\ref{thm:B} does not answer this.
\end{question}

\begingroup\makeatletter\let\addcontentsline\@gobblethree
\section*{Acknowledgments and AI disclosure}
We wish to thank Gil Kalai and Karim Adiprasito for pointing us to valuable references. We used Claude Sonnet 5.5 to generate the examples in Section \ref{sec:examples}, and to carry out some of the computations related to Main Theorems \ref{thm:C}, \ref{thm:D} and \ref{thm:E}. All proofs are still man-made or woman-made.
BB is supported by Simons gift MPS-TSM-00002873; GM is supported by Horizon Europe ERC Grant number: 101045750 / Project acronym: HodgeGeoComb.

\endgroup

\end{document}